\documentclass[10pt,a4paper]{amsart}

\usepackage{graphicx}
\usepackage{psfrag}
\include{psfig}
\newcommand{\ncmd}{\newcommand}
\ncmd{\rencmd}{\renewcommand}
\ncmd{\dspst}{\displaystyle }

\ncmd{\preua}{\noindent \mbox{\bf Proof  }}
\ncmd{\preuda}[1]
{\noindent \mbox{\bf Proof of \ref{#1}   }}

\ncmd{\fin}{\hspace*{\fill} 
\quad\hbox{\hskip 1pt\vrule width 4pt height 6pt
          depth 1.5pt\hskip 1pt} \medskip }

\ncmd{\son}{\mbox{$SO_{0}(1,n-1)\/$ }}
\ncmd{\mink}{\mbox{${\bf M}^{n}\/$ }}

\newtheorem{prop}{Proposition}[section]
\newtheorem{thma}[prop]{Theorem}
\newtheorem{lema}[prop]{Lemma}
\newtheorem{cora}[prop]{Corollary}
\newtheorem{defina}[prop]{Definition}
\newtheorem{raca}[prop]{Remark}
\newtheorem{rica}[prop]{Example}

\newtheorem{banane}{Figure}

\ncmd{\exe}{\begin{ric} \em }
\ncmd{\eexe}{\em \end{ric}}
\ncmd{\rque}{\begin{rac} \em}
\ncmd{\erque}{\em \end{rac}}

\ncmd{\exea}{\begin{rica} \em }
\ncmd{\eexea}{\em \end{rica}}
\ncmd{\rquea}{\begin{raca} \em}
\ncmd{\erquea}{\em \end{raca}}

\newcounter{inc}
\rencmd{\theequation}{\arabic{inc}}

\makeindex

\begin{document}

\title{Globally hyperbolic flat spacetimes}
\author{Thierry Barbot}
\date{February 4, 2004. Work partially supported 
by CNRS and ACI ``Structures g\'eom\'etriques et trous noirs''.}

\begin{abstract}
We consider (flat) Cauchy-complete GH spacetimes, i.e., 
globally hyperbolic flat lorentzian manifolds admitting some
Cauchy hypersurface on which the ambient lorentzian metric restricts
as a complete riemannian metric. We define a family of
such spacetimes - model spacetimes - including four subfamilies: 
translation spacetimes, Misner spacetimes,
unipotent spacetimes, and Cauchy-hyperbolic spacetimes (the
last family - undoubtfull the most interesting one - is a
generalization of standart spacetimes defined by G. Mess).
We prove that, up to finite coverings and (twisted) products by
euclidean linear spaces, any Cauchy-complete GH spacetime 
can be isometrically embedded in a model spacetime, or
in a twisted product of a Cauchy-hyperbolic spacetime by
flat euclidean torus. We obtain as a corollary the classification
of maximal GH spacetimes admitting closed Cauchy hypersurfaces.
We also establish the existence of CMC foliations on every
model spacetime. 
\end{abstract}

\maketitle

\section{Introduction}
Results concerning flat lorentzian manifolds is the mathematical litterature
are mostly devoted to the case of \emph{closed} manifolds (i.e.
compact without boundary). And as a byproduct of this activity, 
the structure of closed flat manifolds has been elucidated in a quite
satisfying way: see for example \cite{barzeg} for a quick survey
on this problem.

But the non-compact case remains essentially open and non considered,
with the notable exception of the works related to
Margulis spacetimes emerging from a question by Milnor
concerning the (non)solvability of discrete groups acting
on Minkowski space (see also \cite{barzeg}).

On the other hand, there is an important and natural notion 
in lorentzian geometry: 
the \emph{global hyperbolicity}.
This notion is central in physics, in the area of classical General
Relativity. And it is incompatible with compactness:
a globally hyperbolic lorentzian manifold is never compact.
Because of this physical background, besides the notational convenience
it provides, we call (flat) lorentzian manifolds
(flat) spacetimes.

A globally hyperbolic (abbrevation GH) spacetime
is a spacetime $M$ admitting a Cauchy hypersurface,
i.e. a hypersurface $S$ which:

- is spacelike (i.e., the lorentzian metric restricts on $S$
as a riemannian metric),

- disconnects $M$,

- intersects every unextendible nonspacelike curve (i.e.,
for which every tangent vector has nonpositive norm)
(see \S \ref{ghdef} for a more complete definition).

$(M,S)$ is a \emph{maximal} GH spacetime (abbrevation MGH)
if the only GH spacetime containing $M$ and for which 
$S$ is still a Cauchy hypersurface, is $M$ itself.
A fundamental theorem by Choquet-Bruhat-Geroch
(\cite{chbr}) states than every GH pair $(M,S)$ 
can be extended in a unique way to
a maximal GH spacetime $N$ in which $S$ is still a
Cauchy hypersurface. 

It is quite surprising that the classification of
flat MGH spacetimes, which is the central topic here, 
has not been previously
systematically undertaken. Such a classification
has its physical interest, and even more, a mathematical one.
It appears as a general extension of Bieberbach's theory to the
lorentzian context.

Actually, such a classification is not possible without some additional
requirement: we will only consider here Cauchy-complete
GH spacetimes, i.e., GH spacetimes admitting a Cauchy hypersurface
on which the lorentzian ambient metric restricts as 
a \emph{complete} riemannian metric. Among Cauchy-complete
GH spacetimes, we distinguish the important subfamily of
Cauchy-compact ones, for which the Cauchy hypersurface
is closed. Cauchy-complete GH spacetimes
which cannot be isometrically embedded in any bigger
Cauchy-complete GH spacetimes, with no additionnal
restriction on Cauchy surfaces,
is \emph{absolutely maximal} (abbrevation: AMGH spacetimes).

The most obvious examples of Cauchy-complete MGH spacetimes
are simply quotients of \mink itself by discrete groups 
of spacelike translations. We call these examples
\emph{translation spacetimes.}

Next, there is another natural family that
we call \emph{Misner spacetimes} (see \S \ref{misner}):
there are quotients of the future in \mink of a spacelike
$(n-2)$-subspace $P$ by an abelian discrete group whose elements
all admit
as linear part a boost (or loxodromic element, maybe trivial) 
acting trivially in $P$, and translation part in $P$.

There is also the family of \emph{unipotent spacetimes,\/}
described in \S \ref{unipfam}: let's briefly mention
here that each of these spacetimes is the quotient of domain
$\Omega$ in \mink delimited by one or two parallel degenerate
hyperplanes by an abelian discrete group of unipotent
isometries of \mink. 

Finally, and maybe the most interesting one, we have the
family of what we call here \emph{Cauchy-hyperbolic spacetimes.\/}
In this family of examples, the linear part 
$L: \Gamma \rightarrow SO_{0}(1,n-1)$
of the holonomy morphism will be injective, with image a discrete subgroup 
of $SO(1,n-1)$ without torsion. Therefore,
we consider $\Gamma \approx L(\Gamma)$ directly as a discrete subgroup
of $SO(1,n-1)$. We moreover assume $\Gamma$ nonelementary, meaning
that the orbits of its action on $\overline{\mathbb H}^{n-1}$ are
all infinite.
 
We first consider the case where $\Gamma$ is a cocompact
lattice of $SO_{0}(1,n-1)$.
Let $\Omega^{+}$ be the connected component of 
$\{ Q < 0 \}$ geodesically complete in the future, and
$\Omega^{-}$ the other connected component: $\Omega^{-}$ is geodesically
complete in the past.
The action of $\Gamma$ on $\Omega^{\pm}$ 
is free and properly discontinuous: we denote
by $M^{\pm}(\Gamma)$ the quotient manifold. Every level set
$\{ Q= -t^{2} \} \cap \Omega^{\pm}$ is $\Gamma$-invariant; it induces in $M^{\pm}(\Gamma)$
a hypersurface with induced metric of constant sectional curvature $-\frac{1}{t^{2}}$.
Recall that $\{ Q=-1 \}$ is the usual representant of the hyperbolic space,
therefore, the flat lorentzian metric on $M^{\pm}(\Gamma)$ admits the warped product
form $-dt^{2}+t^{2}g_{0}$, where $g_{0}$ is the hyperbolic metric on 
$\Gamma\backslash{\mathbb H}^{n-1}$. We call these examples 
\emph{radiant standart spacetimes.\/}
Observe that $M^{+}(\Gamma)$ (resp. $M^{-}(\Gamma)$)
is geodesically complete in the future (resp. in the past), and that
there is a time reversing isometry between them.

In \cite{mess} or \cite{ande}, it is shown that 
any representation of a $\Gamma$ in $\mbox{Isom}({\bf M}^{n})$
admitting as linear part an embedding onto a cocompact lattice of $SO(1,n-1)$
preserves some future complete - and also a past complete - convex
domain of \mink, in such a manner that the quotients
of these domains by $\rho(\Gamma)$ are
Cauchy-compact AMGH spacetimes
. These spacetimes 
are called by G. Mess and L. Andersson
\emph{standart spacetimes}. 

When the linear part $\Gamma$ is still a nonelementary discrete subgroup, but not
cocompact in $SO(1,n-1)$, the question is slightly more
delicate. In \S \ref{listnonele}, we extend the family
of standart spacetimes to a more general one: the family of
Cauchy-hyperbolic spacetimes (definition \ref{definitioncauchyhyp}).
Briefly speaking, they are still quotients of semicomplete convex
domains of \mink by discrete groups of isometries admitting as
linear part a discrete subgroup. But, even in the radiant case,
i.e. when $\rho(\Gamma)$ preserves a point in \mink, let's say, the
origin, the associated Cauchy-hyperbolic spacetime is not always the quotient
$\rho(\Gamma)\backslash\Omega^{\pm}$ as above, but most ofently some bigger
spacetime.

As a last comment on these examples, in \S \ref{suli} we prove
that when $\Gamma$ is a \emph{convex cocompact Kleinian group}, then
any discrete subgroup of $\mbox{Isom}({\bf M}^{n})$ admitting
$\Gamma$ as linear part preserves a semicomplete domain of
\mink, i.e. is the holonomy group of some Cauchy-hyperbolic
spacetime. This claim is not as trivial as it may appear
at first glance: in proposition \ref{tri}, we exhibit
a (nonuniform) lattice $\Gamma$ in $SO_{0}(1,2)$) 
such that, a representation
$\rho: \Gamma \rightarrow \mbox{Isom}({\bf M}^{n})$
for which $L \circ \rho$ is the identity morphism can preserve
a semicomplete convex domain if and only if it is radiant,
i.e. preserves a point in \mink. This fuchsian group 
$\Gamma$ is nothing but the group associated to
the $3$-punctured conformal sphere.

Convex cocompact Kleinian groups form an important family
including cocompact Kleinian groups or Schottky groups, which are essentially
the only examples appearing in the physical litterature, except geometrically
finite Kleinian groups (convex cocompact Kleinian groups can be defined
as geometrically finite Kleinian groups without parabolic elements). 
The correct way to extend this result to geometrically finite Kleinian 
groups is an interesting question,
even in the $2+1$-dimensional case.

Before stating the classification's Theorems, let's indicate
a natural way to produce GH spacetimes from other
ones: given a flat GH spacetime $M$,
a flat euclidian manifold $N$, 
and a representation $r: \pi_{1}(M) \rightarrow Isom(N)$, 
the total space $B$
of the suspension bundle over $M$ 
with fiber $N$ and monodromy $r$, is naturally equipped
with a flat lorentzian metric which is still
globally hyperbolic. The resulting spacetime $B$ from this construction 
which is defined with more details in \S \ref{produi}, is called the 
twisted product of $M$ by $N$ with monodromy $r$.
The case where $N$ is an euclidean linear space and
the representation $r$ admits a global fixed point in $N$
is particularly pleasant: we say then that $B$ is a linear
twisted product over $M$.

We now express the main theorem of this paper, essentially
stating that the examples above provide 
the complete list of flat Cauchy-complete 
MGH spacetimes. In the following statement, a \emph{tame embedding}
is an isometric embedding inducing an isomorphism between
fundamental groups.

\begin{thma}
\label{caucomp} 
Up to finite coverings and linear twisted products, every Cauchy-complete
globally hyperbolic flat spacetime can be tamely embedded
in an absolutely maximal globally hyperbolic spacetime 
which is a translation spacetime, a
Misner spacetime, a unipotent spacetime,
or the twisted product of a  Cauchy-hyperbolic flat GH spacetime
by an euclidean torus.
\end{thma}

The Cauchy-compact case deserves its own statement:

\begin{thma}
\label{cauclo}
Up to finite coverings,
every maximal Cauchy-compact globally hyperbolic flat spacetime is isometric
to a translation spacetime, a Misner spacetime, or 
the twisted product of a standart spacetime (Cauchy-compact 
Cauchy-hyperbolic GH spacetime) by an euclidean torus.
\end{thma}

The proofs of these theorems are written here in a quite intricate
way, since it involves particular subcases to consider separately. Hence,
we collect along the text all the elements of the proof, and
then indicate in \S \ref{sumup} the way to reconstruct from these
intermediate results the complete proofs of the theorems.

Theorem \ref{caucomp} actually is not a full classification theorem: indeed,
in a given MGH spacetime we can embed many different
MGH spacetimes - consider for example the Minkowski space itself!
But, such a complete classification of Cauchy-complete
GH spacetimes is not only unrealistic, it is furthermore 
useless when we adopt the point of view that the object
of study here is the holonomy group $\Gamma$.
In other words, our essential procedure is to
associate to suitable discrete subgroups of $\mbox{Isom}({\bf M}^{n})$ 
some invariant domains of \mink (regular convex domains) on which
the dynamical behaviour of the group respect some causality
properties.

As a corollary of theorem \ref{cauclo}, 
we obtain that closed Cauchy hypersurfaces of 
globally hyperbolic flat spacetimes are homeomorphic to
finite quotients of products 
of torii with hyperbolic manifolds.
This generalizes \cite{scannel} in any dimension, without the superfluous 
$3$-dimensional topological arguments used in \cite{scannel}.

There are other works related to the present work: 
in \cite{scannel2}, K. Scannell classified 
Cauchy-compact maximal globally hyperbolic spacetimes
with constant sectionnal curvature $+1$:
the topological type of the Cauchy hypersurface $S$ being fixed, there is a $1-1$ correspondance
between MGH spacetimes and riemannian flat conformal structures on $S$. 
Maximal globally hyperbolic Cauchy-compact spacetimes with constant curvature $-1$ of
dimension $2+1$ are classified in \cite{mess}. We should mention that results in \cite{scannel},
\cite{mess} are not stated in the terminology of GH spacetimes, but our presentation
follows immediatly from their works.

In \cite{ande}, L. Andersson
classified flat Cauchy-compact MGH spacetimes, but with the initial hypothesis
that Cauchy hypersurfaces have hyperbolic type. He proves also that these spacetimes
all admit foliations by constant mean curvature hypersurfaces (abbreviation,
CMC foliation), which is unique in a given spacetime. Thanks to theorem \ref{cauclo},
we can extend this result to the elementary case. It is done in \S \ref{CMC}.

The present paper include also a generalization to any dimension
of the characterization of "spacelike regions" for isometries of
\mink made in \cite{drumm} 
("spacelike regions" of \cite{drumm} are
called here achronal domains).

\section{Preliminaries}

\subsection{Spacetimes}
By spacetime, we mean here in general a (non closed) oriented
and chronologically oriented lorentzian manifold. 
Our convention here is that lorentzian metrics have
signature $(-, +, \ldots, +)$. Tangent vectors are called
\emph{spacelike}, \emph{timelike}, \emph{lightlike} if their norms
are respectively positive, negative, null. 
A curve in the spacetime is nonspacelike if its tangent 
vectors have nonpositive norm.
In the same spirit, a curve with 
timelike tangent vectors is called timelike, it 
has to be thought as a potential trajectory
of a particle. Its proper time is defined by:

\[ \mbox{proper-time}(c) = \int \sqrt{-\langle \partial_{t}c \mid \partial_{t}c
\rangle} \]

Since the spacetime is chronologically oriented, every 
nonspacelike curve admits a canonical orientation towards
its future.

A spacetime is geodesically complete if every
timelike geodesics admit geodesic parame\-trizations by
$]-\infty, +\infty[$. It is \emph{geodesically complete in the future}
(resp., in the past) if every future oriented (resp. past oriented)
timelike geodesic ray admits
a geodesic parametrization by $]0, +\infty[$.
Finally, a \emph{geodesically
semicomplete\/} spacetime is a spacetime which is either geodesically
complete in the future or geodesically complete in the past.

\subsection{Minkowski space}
We denote here by \mink the Minkowski space.
We stress out that we consider here \mink endowed
with an orientation and a chronological orientation.

The lorentzian quadratic form on \mink is denoted by 
$\langle . \mid . \rangle$; and the minkowski norm is 
$\mid x \mid^{2} = \langle x \mid x \rangle$ (the notation
$\mid x \mid$ will be reserved to spacelike vectors).
For any affine subspace $E$ of \mink, we denote by
$E^{\perp}$ the orthogonal of $E$ - this is a subspace
of the underlying linear space. When $E$
is a lightlike affine line, $E^{\perp}$ will also
abusively denote the unique degenerate affine hyperplane
containing $E$ with direction orthogonal to $E$.

Let $L: \mbox{Aff}(n, {\mathbb R}) \rightarrow GL(n, {\mathbb R})$ be the 
usual linear part morphism for the group of affine
transformations of \mink. The isometry group of the 
(unoriented) Minkowski space is the space of affine
transformations with linear part in $O(1,n-1)$, but
since here the consider \mink as 
oriented and chronologically orientated, the term
isometry will be reserved to affine transformations 
admitting linear parts in the identity
component $SO_{0}(1,n-1)$ (also called the orthochronous component). 
We denote by $\mbox{Isom}({\bf M}^{n})$ the group
of isometries of \mink.

\subsection{Flat spacetimes as geometric manifolds}
Truely speaking, every spacetimes considered in this work
are flat, i.e., locally modeled on the Minkowski space.
In other words, in the language of geometric structures
(see e.g. \cite{gold}),
they are $(G,X)$-manifolds with $X = {\bf M}^{n}$ and
$G = \mbox{Isom}(M^{n})$.

We will need only few facts on geometric structures: 
the existence of the developping map and the holonomy
morphism (see \S \ref{spapa}).

\subsection{Globally hyperbolic spacetimes}
\label{ghdef}
A spacetime $M$ is \emph{globally hyperbolic\/} 
(abbreviation GH) if there is a proper time function
$t: M \rightarrow \mathbb R$  
such that every fiber $t^{-1}(t_{0})$
is spacelike. Moreover, it is required that any nonspacelike curve in $M$ can
be extended to another one such that the restriction of $t$ on this extended
curve is a diffeomorphism onto the entire $t(M)$. In other words,
every fiber of $t$ is a \emph{Cauchy hypersurface,\/} meaning precisely that
this is a hypersurface intersecting 
every unextendible nonspacelike curve in $M$.
There are other equivalent definitions of global hyperbolicity, see \cite{beem}.
Observe that $t: M \rightarrow {\mathbb R}$ is necessarly a locally trivial, thus
trivial fibration: $M$ is diffeomorphic to $S \times {\mathbb R}$.
Observe also that GH spacetimes admit chronological orientations.

Choquet-Bruhat-Geroch's work (\cite{chbr}) implies that every
GH manifold admits
a unique maximal globally hyperbolic extension. Let's be more precise:
let $S$ be a Cauchy hypersurface in a globally hyperbolic spacetime $M$.
A $S$-embedding of $M$ is a isometric embedding $f: M \rightarrow M'$
where $M'$ is a spacetime, such that $f(S)$ is still a Cauchy hypersurface
in $M'$. Actually, this notion does not depend on the choice of
the Cauchy hypersurface: if $S$, $S'$ are Cauchy hypersurfaces in $M$,
an isometric embedding $f: M \rightarrow M'$ is a $S$-embedding if and
only if it is a $S'$-embedding. Therefore, we call such a map
a \emph{Cauchy embedding.\/} Now, a globally hyperbolic manifold
$M$ is said maximal if any Cauchy embedding of $M$ into some globally hyperbolic
spacetime $N$ is necessarly surjective. Choquet-Bruhat-Geroch's result can now
be precisely stated: every globally hyperbolic spacetime $M$ admits
a Cauchy embedding in a maximal globally hyperbolic spacetime. Moreover,
this maximal globally hyperbolic extension is unique up to isometries.

\rquea
The maximal globally hyperbolic extension of a flat globally hyperbolic spacetime
is flat. Observe also that if the spacetime is analytic, with analytic
lorentzian metric, then it is flat as soon it contains a flat open set.
\erquea

\rquea MGH spacetimes may
admit nonsurjective isometric embeddings in 
bigger MGH spacetimes. Consider the following example:
in \mink, let $\Omega$ be a connected component of the negative cone $\{ Q< 0 \}$.
Then, $S = \Omega \cap \{ Q=-1 \}$ is a Cauchy hypersurface for $\Omega$.
Observe that $\Omega$ is MGH, but nontrivially embedded in ${\bf M}^{n}$!
The point is that $S$ is not a Cauchy hypersurface in \mink, therefore,
the inclusion $\Omega \subset {\bf M}^{n}$ is not a Cauchy embedding.
\erquea

\subsection{Cauchy-complete, Cauchy-compact GH spacetimes}
A GH spacetime is \emph{Cauchy-complete} if it admits a Cauchy hypersurface
where the ambiant lorentzian metric restricts as a complete riemanian
metric. This notion depends on the choice of 
the Cauchy hypersurface: for example, in the two-dimensional Minkowski
space, i.e. the plane equipped with the metric $dxdy$, Cauchy hypersurfaces
are graphs $y = f(x)$ of strictly increasing diffeomorphisms from $\mathbb R$ onto
$\mathbb R$, and such a graph is complete if and only if the integrals of the square root
of $f'$ on $]-\infty, 0]$ and $[0, +\infty[$ are infinite.
 
Among Cauchy-complete spacetimes, there is a remarkable natural
family: GH spacetimes with \emph{closed\/} Cauchy hypersurfaces (we know
that in general, if a spacetime admits a closed Cauchy hypersurface, then
all its Cauchy surfaces are closed - indeed, they are homeomorphic one to
the other). We call such a spacetime \emph{Cauchy-compact\/} GH spacetime
(abbreviation CGH). It is worth to emphasize this family, since our results
are much easier to state in this situation; keeping in mind our
analogy with Bieberbach's theory on complete euclidean manifolds,
Cauchy-compact spacetimes correspond to closed euclidean manifolds.

\rquea
\label{clomax}
We say that a flat Cauchy-complete GH spacetime $M$ is absolutely
maximal (abbreviation AMGH) if any isometric embedding into another 
Cauchy-complete flat GH spacetime is surjective. 
Cauchy-compact GH spacetimes are
automatically absolutely maximal. Indeed, any spacelike 
hypersurface in a given CGH spacetime is a Cauchy hypersurface.
\erquea

\section{Elementary Cauchy-complete AMGH spacetimes}
\label{listele}

\subsection{Translation spacetimes}
These are the easiest and more obvious examples, including
\mink itself: they are the
quotients of the entire \mink by a group of translations by spacelike
translation vectors. 
There are orthogonal sums $S \oplus {\mathbb R}$, where $S$ is a flat euclidean
cylinder ${\mathbb T}^{k} \times {\mathbb R}^{n-k-1}$ (where ${\mathbb T}^{k}$ is a flat
torus of dimension $k$), and ${\mathbb R}$ is equipped with the metric $-dt^{2}$.

\subsection{Misner flat spacetimes}
\label{misner}
Consider the metric $2dxdy+dz^{2}$ on ${\mathbb R}^{n}$, where $z$ denotes an element
of ${\mathbb R}^{n-2}$, and an isometry $\gamma_{0}$
$(x,y,z) \mapsto (e^{t_{0}}x, e^{-t_{0}}y, z ) \;\;\; (t_{0}\neq 0)$. The
$1$-dimensional Misner spacetime is
the quotient of $\Omega = \{ xy<0 , x>0 \}$ by the group 
generated by $\gamma_{0}$. Consider the pure loxodromic
element $\gamma_{0}$ above as the time $1$ of a $1$-parameter subgroup 
$\gamma_{0}^{t}$, and let $G_{0}$ be the subgroup of $\mbox{Isom}({\bf M}^{n})$
formed by elements with linear parts in $\gamma_{0}^{t}$ and translations
parts in the subspace $\{ x = y = 0 \}$. This an abelian Lie group
isomorphic to ${\mathbb R}^{n-1}$, and the orbits of $G_{0}$ in $\Omega$
are all spacelike. If $\Gamma$ is any discrete subgroup of
$G_{0}$, the quotient $\Gamma\backslash\Omega$ is a flat globally hyperbolic
spacetime, that we call \emph{Misner spacetime.\/}

Observe that the linear part $L(\Gamma)$ is not necessarely discrete:
it can be a dense subgroup of $\gamma_{0}^{t}$. Observe also that
we could have selected $\Omega = \{ xy<0 , x<0 \}$ too; but the two choices
provide isometric spacetimes (up to a time reversing isometry), which 
are both geodesically semicomplete.

As a last comment on this example, we should add a particularly comfortable
coordinate system: paramatrize $\Omega$ by $x=e^{\eta+\nu}$, $y=-e^{\eta-\nu}$.
The flat metric $2dxdy+dz^{2}$ on $\Omega$ in these coordinates is:

\[ e^{2\eta}( -2d\eta^{2}+2d\nu^{2}+ e^{-2\eta}dz^{2}) \]

The action of $G_{0}$ in this coordinate system is the action by translations
on the $z, \nu$ coordinates.

\subsection{Unipotent spacetimes}
\label{unipfam}
Consider once more a coordinate system $(x,y,z)$ with $x$, $y$ in ${\mathbb R}$
and $z$ in ${\mathbb R}^{n-2}$, where the metric is given by $2dxdy+dz^{2}$.
We consider the unipotent part of the stabilizer of the lighlike hyperplane
$\{ y=0 \}$. This is the group $\mathcal A$ with elements of the form:

\[ g_{u,v,\mu}(x,y,z) = (x+\mu-\langle z \mid v \rangle - \frac{1}{2}y\mid v \mid^{2},
y, z+u+yv) \]

where $u$ and $v$ are elements of ${\mathbb R}^{n-2}$, and $\mu$ a real number.
The orbits of $\mathcal A$ are the degenerate hyperplanes
$y=Cte$. Observe that $\mathcal A$ is not abelian, but a central extension
of $\mathbb R$ by ${\mathbb R}^{2n-4}$: the operation law is:

\[ g_{u,v,\mu} \circ g_{u',v',\mu'} = g_{u+u',v+v',\mu+\mu'-\langle u' \mid v \rangle} \]

Select now $n-2$ real numbers $\lambda_{1} \leq \lambda_{2} \leq \ldots \leq \lambda_{n-2}$,
and an orthogonal basis $e_{1}, e_{2}, \ldots, e_{n-2}$ of the euclidean space
${\mathbb R}^{n-2}$. Then, when $(t_{1}, \ldots , t_{n-2})$
describe the entire ${\mathbb R}^{n-2}$, the $g_{u,v, \mu}$ with $u = \sum t_{i}e_{i}$, 
$v = \sum t_{i}\lambda_{i}e_{i}$ and 
$\mu = -\sum t_{i}^{2}\frac{\lambda_{i}}{2}$ describe an abelian subgroup $A$
of $\mathcal A$ isomorphic to ${\mathbb R}^{n-2}$. 

The action is given by:

\[ (t_{1}, \ldots, t_{n-2}).(x,y,\sum z_{i}e_{i}) = (x-\sum \frac{{\lambda}_{i}}{2}t_{i}^{2}-\sum \lambda_{i}z_{i}t_{i}-\frac{y}{2}\sum \lambda_{i}^{2}t_{i}^{2}, y, \sum (z_{i}+t_{i}(1+y\lambda_{i}))e_{i}) \]

The domain of
\mink where the orbits of $A$ are spacelike
is $\Omega(A) = \{ (x,y,z) \slash 1+y\lambda_{i} \neq 0 \;\; \forall i  \}$.

Define $y_{i} = -\frac{1}{\lambda_{i}}$ ($y_{i}=\infty$ if $\lambda_{i}=0$),
and  
$y_{0}=-\infty \leq y_{1} \leq \ldots y_{n-2} \leq y_{n-1} =+\infty$. 
Then, a connected component $\Omega$ of 
$\Omega(A)$ is the set of $(x,y,z)$ with
$y_{j} < y < y_{j+1}$ for some index $j$. 

Consider such a connected component, and a discrete
subgroup $\Gamma$ of $A$. Observe that the boundary of $\Omega$ is the union of one or 
two lightlike hyperplanes $\{ y=Cte \}$. 

Let $f: I=]y_{j}, y_{j+1}[ \rightarrow {\mathbb R}$ be any 
(strictly) increasing function. When $y$ describes $I$,
$(f(y), y, 0)$ describes a spacelike curve $L$ in $\Omega$. The saturation of $L$ under
the action of $A$ is the set of $(x,y,z)$ with:

- $z=\sum t_{i}(1+\lambda_{i}y)e_{i}$,

- $x=f(y)-\sum \frac{{\lambda}_{i}}{2}t_{i}^{2}-\frac{y}{2}\sum \lambda_{i}^{2}t_{i}^{2}$

Thus, it is the graph of the function
$\Phi(y,z) = f(y)+\frac{1}{2}\sum \frac{z_{i}^{2}}{y_{i}-y}$ 
(with the conventions $z=\sum z_{i}e_{i}$, and 
$\frac{z_{i}^{2}}{y_{i}-y}=0$ if $y_{i}=\infty$).
Of course, $S_{f}$ is a spacelike $A$-invariant hypersurface. More precisely,
the Minkowski norm $2dxdy+\sum dz_{i}^{2}$ induces on $S_{f}$ the metric:

\[ 2f'(y) dy^{2} + \sum (dz_{i}-\frac{z_{i}dy}{y-y_{i}})^{2} \]

with the convention $\frac{z_{i}dy}{y-y_{i}} = 0$ if
$y_{i} = 0$.
If $\zeta_{i}=\frac{z_{i}}{y-y_{i}}$ (convention: $\zeta_{i}=z_{i}$ if $y_{i} = \infty$), the expression in the 
$(y, \zeta_{1}, \ldots, \zeta_{n-2})$ coordinates of the metric on $S$ is
(with the convention $(y-y_{i})d\zeta_{i} = d\zeta_{i}$ if $y_{i} = \infty$):

\[ 2f'(y) dy^{2} + \sum (y-y_{i})^{2}d\zeta_{i}^{2} \]

This is time now to specify the function $f$: let $y_{0}$ be any 
element of $I$, and $a: I \rightarrow ]1, +\infty[$ any 
smooth function such that the integrals $\int_{y_{j}}^{y_{0}} a(y) dy$ 
and $\int_{y_{0}}^{y_{j+1}} a(y) dy$ are both infinite. 
Choose $f(y) = \int_{y_{0}}^{y}a^{2}(y)dy$. We claim that for such a $f$, 
the metric on $S_{f}$ is complete:
if not, there is a smooth path $c: [0, +\infty[ \rightarrow S_{f}$ 
with finite length and escaping to infinity.
We write $c(t)=(y(t), \zeta_{1}(t), \ldots, \zeta_{n-2})$.
Since the length is finite, the integrals 
$\int_{0}^{+\infty} \sqrt{2f'(y(t))}\mid y'(t)\mid dt$ and 
$\int_{0}^{+\infty} \sqrt{\sum (y(t)-y_{i})^{2}\zeta'_{i}(t)^{2}} dt$ 
are both finite. The finiteness of the first
integral implies the finiteness of $\int_{0}^{+\infty}a(y(t))y'(t)dt$. 
It follows that the $y$-coordinate remains in a compact subintervall of $I$. 
But, then, the integral $\int_{0}^{+\infty} \sqrt{\sum \zeta'_{i}(t)^{2}} dt$ 
is bounded. It follows that $c(t)$ stay in a compact domain
of $S_{f}$: contradiction.

Hence, $S_{f}$ is a complete spacelike hypersurface. 

We now prove that $S_{f}$ is a Cauchy hypersurface for $\Omega$. 
Clearly, the complement of $S_{f}$ in \mink has two connected 
components, one containing $\{ y=y_{j} \}$, the other
containing $\{ y=y_{j+1} \}$. Now, a nonspacelike geodesic 
in $\Omega$ is the intersection between $\Omega$ and a 
nonspacelike geodesic line $d$ of \mink. There are two cases to consider:
if the direction of $d$ is the isotropic direction $\Delta_{0}$,
then $d$ must intersect $S_{f}$ since $S_{f}$ is a graph $x=\Phi(y,z)$.
If not, $d$ intersect the two lightlike hyperplanes 
$\{ y=y_{j} \}$,  $\{ y=y_{j+1} \}$.
Then, it must intersect $S_{f}$ since $S_{f}$ disconnects these hyperplanes.

Hence, any nonspacelike geodesic intersect $S_{f}$: this is a criteria proving that 
$S_{f}$ is indeed a Cauchy hypersurface.

It follows that the quotient of $\Omega$ by $\Gamma$ is globally
hyperbolic and Cauchy-complete.
We call such a spacetime \emph{unipotent spacetime.\/} 
Unipotent spacetimes have \emph{never} closed Cauchy hypersurfaces, 
since the $\mathcal A$-invariant coordinate
$y$ defines a submersion from
the Cauchy hypersurface into ${\mathbb R}$. Observe also that unipotent spacetimes are
not uniquely defined by their holonomy: indeed, the same $\Gamma$ may be extended
to a $A \approx {\mathbb R}^{n-2}$ in different ways. Moreover, even when $A$
is fixed, different choices of connected components of $\Omega(A)$ lead
to nonisometric unipotent spacetimes. Finally, unipotent
spacetimes are not geodesically semicomplete, except whose corresponding
to connected components $\Omega$ with only one boundary hyperplane.

\section{The non-elementary GH spacetimes: Cauchy-hyperbolic spacetimes}
\label{listnonele}

\subsection{The Penrose boundary}
\label{penrose}
In this section, we describe ${\mathcal J}$, 
the space of lightlike (or degenerate)
affine hyperplanes of \mink. 
We represent as usual \mink as ${\mathbb R}^{n}$
equipped with the metric $-dx_{0}^{2}+dx_{1}^{2} +  \ldots + dx_{n-1}^{2}$.
We will also use the euclidean metric 
$dx_{0}^{2}+dx_{1}^{2} +  \ldots + dx_{n-1}^{2}$; the norm of
a vector $v$ of ${\mathbb R}^{n}$ for this metric will be denoted $N(v)$.

Denote by $\mathcal S$ the space of linear lightlike hyperplanes.
It admis a natural
$1-1$-parametrisation by vectors $v$ of 
${\mathbb R}^{n}$ for which $N(v)=1$, 
$\mid v \mid = 0$, and $v$ oriented towards the future,
i.e., $x_{0} > 0$. We denote by ${\mathcal S}^{+}$
the space of such vectors. A vector in ${\mathcal S}^{+}$ is 
completely characterized by
its $x_{1}, \ldots, x_{n-1}$ coordinates which satisfy 
$\sum x_{i}^{2} = \frac{1}{2}$ (observe that necessarely
$x_{0} = \frac{1}{2})$). Hence, ${\mathcal S}^{+}$ is naturally
identified with the $(n-2)$-sphere. 
The ambient lorentzian metric on \mink
induces on $\{ N=1 , \mid . \mid = 0 \}$ a non-degenerate
riemanian metric, and an easy computation proves that this
metric is nothing but the usual round metric on the
$(n-2)$-sphere. In the same vein, ${\mathcal S}$ can
also parametrized by ${\mathcal S}^{-}$, the space of
\emph{past oriented} lightlike vectors with $N$-norm $1$.

The group $SO_{0}(1, n-1)$ acts naturally on ${\mathcal S}$
but it does not preserve the ``round metric'' we have just defined -
mainly because it does not preserves $N=1$. 
Actually, for $g$ in $SO_{0}(1,n-1)$, and identifying
${\mathcal S}^{\pm}$ with the space of $v$'s as above, the
action of $g$ on $\mathcal S$ maps $v$ on $[g].v = \frac{g(v)}{N(g(v))}$.
This action preserves the \emph{conformal class} of the metric;
in particular, it preserves the angles. As a matter of
fact, ${\mathcal S}^{\pm}$ equipped with this conformal class
is the usual conformal sphere, and elements of
$SO_{0}(1,n-1)$ act on it as M\"{o}bius transformations.

The space ${\mathcal J}$ admits a natural fibration 
over each ${\mathcal S}^{\pm}$: this map is just the map associating
to the affine hyperplane its direction. We denote it
by $\delta^{\pm}: {\mathcal J} \rightarrow {\mathcal S}^{\pm}$.
Now, for any element $P$ of ${\mathcal J}$ with
$v = \delta^{\pm}(P)$, the scalar product $\langle v \mid p \rangle$
does not depend on the point $p$ of $P$: we denote it by
$\nu^{\pm}(P)$. Then, the product map 
$\delta^{\pm} \times \nu^{\pm}: {\mathcal J} \rightarrow {\mathcal S}^{\pm} \times
{\mathbb R}$ is one-to-one. 

\rquea
\label{ambigu}
The $\pm$-ambiguity means that we actually have two natural identifications
of $\mathcal J$ with the product of the conformal sphere by the real
line. We will indicate which of these identifications map 
$\delta^{\pm} \times \nu^{\pm}$ is considered by denoting
the space of lightlike hyperplanes by ${\mathcal J}^{+}$ or
${\mathcal J}^{-}$. Alternatively, we can adopt the point
of view that ${\mathcal J}^{+}$ (resp. ${\mathcal J}^{-}$)
is the space of future complete (resp. past complete) half
spaces bounded by lighlike hyperplanes.
\erquea

Finally, $Isom({\bf M}^{n})$ acts naturally on ${\mathcal J}$.
When we use the identification 
${\mathcal J} \approx {\mathcal J}^{+} \approx {\mathcal S} \times {\mathbb R}$, the action
of an isometry $g(x) = L(g)x + \tau$ is expressed by:

\[ g.( v , s ) = ( [L(g)].v , a(g, v)s + b(g, v) ) \]

where:

\begin{eqnarray*}
[L(g)].v & = & \frac{L(g)v}{N(L(g)v)} \\
a(g, v) & = & \frac{1}{N(L(g)(v))} \\
b(g, v) & = & \frac{ \langle \tau \mid L(g)v
\rangle}{N(L(g)v)}
\end{eqnarray*}

Observe that for the other identification 
${\mathcal J} \approx {\mathcal J}^{-}$, the action would be expressed
in a similar form, with the same $a(g, v)$, but with an oppopsite
$b(g, v)$. Observe also that the factor 
$\frac{1}{N(L(g)v)}$ is common to the $v$ and
$s$ components. Hence:

\begin{prop}
\label{laoucapulse}
Equip ${\mathcal J}^{+} \approx {\mathcal S}^{+} \times {\mathbb R}$
with the product metric of the round metric on the sphere ${\mathcal S}^{+}$
and any euclidean metric on $\mathbb R$. Let $g$ be a linear
isometry of \mink, i.e., admitting a trivial translation part
$\tau$. Then, 
the domain in ${\mathcal J}^{+}$
where the action of $g$ is expanding is precisely the preimage under
$\delta^{+}$ of the domain of ${\mathcal S}^{+}$ where the action of $L(g)$
is expanding.
\end{prop}

\rquea
The title of this section is justified by the fact that
${\mathcal J}^{\pm}$ are nothing but the regular parts 
of the Penrose boundary of Minkowski space as usually defined
in general relativity (and usually with the same notation, 
see e.g. \cite{ellis}).
\erquea

\subsection{Regular convex domains}

Let $\Lambda$ be any subset of ${\mathcal J}$.
We will always assume that $\Lambda$ contains at least one
point.
For any element $P$ of $\Lambda$, let $P^{+}$
be the future of $P$, and $P^{-}$ the past
of $P$. These are half-spaces, admitting both $P$
as boundaries. If $v = \delta^{+}(P)$ and $s = \nu^{+}(P)$, 
$P^{+}$ (resp. $P^{-}$) is the domain of points $p$
in \mink for which $\langle v \mid p \rangle - s$ is negative
(resp. positive). 

\begin{defina}
The future complete (resp. past complete) 
convex set defined by $\Lambda$ is the intersection:

\[ \Omega^{+}(\Lambda) = \bigcap_{P \in \Lambda} P^{+} \;\;\;\;\; 
(\mbox{resp.} \;
\Omega^{-}(\Lambda) = \bigcap_{P \in \Lambda} P^{-}) \]

If $\Lambda$ contains at least $2$ elements, then $\Omega^{\pm}(\Lambda)$,
if nonempty, is regular.
\end{defina}

We first collect some straightforward observations:

- $\Omega^{\pm}(\Lambda)$ is an convex set.

- $\Omega^{+}(\Lambda)$ is (geodesically) complete in the future, and 
$\Omega^{-}(\Lambda)$ is complete in the past. They are disjoint one
from the other.

\begin{defina}
If $\Omega^{+}(\Lambda)$ (resp. $\Omega^{-}(\Lambda)$) is
nonempty and open, then it is called a future (resp. past)
regular convex domain. The set $\Lambda$ is then said
future regular (resp. past regular).
\end{defina}

This definition of regular convex domains is completely equivalent to the
definition given by F. Bonsante (see \cite{bonsante}, definition $4.1$),
but it is essential for our purposes to relate it with closed subsets
of ${\mathcal J}$.

\begin{lema}
\label{fermons}
If $\overline{\Lambda}$ is the closure in 
${\mathcal J}$ of $\Lambda$, then $\Omega^{\pm}(\overline{\Lambda})$
contains the interior of $\Omega^{\pm}(\Lambda)$.
\end{lema}

\preua
Let $x$ be a point in the interior of $\Omega^{+}(\Lambda)$.
Assume that $x$ does not belong to $\Omega^{+}(\overline{\Lambda})$.
Then, for some element $(v,s)$ of $\Omega^{+}(\overline{\Lambda})$,
we have $\langle x \mid v \rangle \geq s$. On the other hand,
for any sequence $(v_{n}, s_{n})$ in $\Lambda$ converging
to $(v,s)$, we have $\langle x \mid v_{n} \rangle < s_{n}$.
At the limit, we obtain: $\langle x \mid v \rangle = s$.
But, since $x$ is in the interior of $\Omega^{+}(\Lambda)$,
there is a point $y$ in the interior of $\Omega^{+}(\Lambda)$ near $x$
for which $\langle y \mid v \rangle > s$. Apply once more
the argument above to $y$: we obtain a contradiction.\fin

\begin{cora}
\label{closeregular}
If $\Lambda$ is future regular, then its closure $\overline{\Lambda}$ is future
regular too, and $\Omega^{+}(\Lambda) = \Omega^{+}(\overline{\Lambda})$.
\end{cora}

\preua
The inclusion $\Omega^{+}(\overline{\Lambda}) \subset \Omega^{+}(\Lambda)$ is
obvious. Hence, the corollary follows immediatly from
lemma \ref{fermons}.\fin

\begin{lema}
\label{closecool}
If $\Lambda$ is closed, then $\Omega^{+}(\Lambda)$ is open
(maybe empty).
\end{lema}

\preua
Assume the existence of some element $x$ belonging to $\Omega^{+}(\Lambda)$,
but not to its interior. Then, there is a sequence of elements
$x_{n}$ in \mink converging to $x$, but not belonging to
$\Omega^{+}(\Lambda)$. It means that for every $n$, 
there is an elements $(v_{n}, s_{n})$ in $\Lambda$
for which $\langle x_{n} \mid v_{n} \rangle \geq s_{n}$.
By compactness of $\mathcal S$, we can assume that the
sequence $v_{n}$ converges to some $v$.
Since $x$ belongs to $\Omega^{+}(\Lambda)$, we have:
$\langle x \mid v_{n} \rangle < s_{n}$. It follows
that $s_{n}$ converges to $\langle x \mid v \rangle$.
But since $\Lambda$ is closed, it must contain 
the limit point $(v, \langle x \mid v \rangle)$.
This is a contradiction since $x$ is assumed to
belong to $\Omega^{+}(\Lambda)$.\fin

\begin{lema}
\label{fermonsbis}
Let $\Lambda$ be a subset of ${\mathcal J}$ with
future regular closure. Then, $\Omega^{+}(\overline{\Lambda})$
is the interior of $\Omega^{+}(\Lambda)$.
\end{lema}

\preua
According to lemma \ref{fermons}, $\Omega^{+}(\overline{\Lambda})$
contains the interior of $\Omega^{+}(\Lambda)$. But it is
obviously contained in $\Omega^{+}(\Lambda)$, and, according to
lemma \ref{closecool}, it is open. The lemma follows. \fin

Thanks to corollary \ref{closeregular} and lemma \ref{closecool}, 
we can define regular convex domains as \emph{nonempty} convex sets associated
to \emph{closed} sets $\Lambda$ in ${\mathcal J}$ not reduced 
to one point (observe that the similar statements for
past regular domains are evidently valid).

\begin{prop}
A closed subset of ${\mathcal J}$ is future regular
(resp. past regular) 
if and only if there is some real number $C$ such that,
for every $(v,s)$ in $\Lambda$, the second component $s$
is less than $C$ (resp. bigger than $C$).
\end{prop}

\preua
Assume that $\Lambda$ is contained 
in a domain $\{ s > C \}$. Then, $\{ s=C \}$ is a family 
of lightlike hyperplanes tangent to the future cone
of some point $x$ satisfying $\langle x \mid v \rangle = C$
for every future oriented lightlike $v$ satisfying
$N(v)=1$. In particular, $\langle x \mid v \rangle < s$
for every $(v,s)$ in $\Lambda$:
$x$ belongs to $\Omega^{+}(\Lambda)$
which therefore is nonempty.

Inversely, if $x$ belongs to $\Omega^{+}(\Lambda)$, 
then, for every $(v,s)$ in $\Lambda$,
we have $\langle x \mid v \rangle < s$.
Since $\mathcal S$ is compact, $\langle x \mid v \rangle$
admits a lower bound independant from $v$.
The proposition follows (the past regular case
is completely similar).\fin

\begin{cora}
\label{futetpast}
$\Lambda$ is future regular and past regular if and only
if it is compact.\fin
\end{cora}

In the following, we state many results proved in \cite{bonsante}
about the \emph{cosmological time of $\Omega^{\pm}(\Lambda)$,}
without any attempt to provide proofs. We just suggest to
the reader to keep in mind the similar situation in euclidean space
when is considered the distance function to some convex set.

For any point $p$ in $\Omega^{+}(\Lambda)$, consider every
future oriented timelike curve starting from a point in 
$\partial\Omega^{+}(\Lambda)$ and ending at $p$. The proper times of
all these curves are uniformely bounded:
let $T(p)$ be the supremum value of these proper times. 
There is a unique $\pi(p)$ in $\partial\Omega^{+}(\Lambda)$
for which $T(p)$ is equal to $-\sqrt{\mid p-\pi(p) \mid^{2}}$; in other words,
the straight segment from $\pi(p)$ to $p$ is the unique timelike
curve linking $\partial\Omega^{+}(\Lambda)$ with proper time 
realizing $T(p)$. The function $T: \Omega^{+}(\Lambda) \rightarrow {\mathbb
R}^{+}_{\ast}$
is a $C^{1}$-convex function (proposition $4.3$ in \cite{bonsante}).
The level sets $S_{t} = T^{-1}(t)$ are convex spacelike hypersurfaces (corollary
$4.5$ in \cite{bonsante}); actually, the direction of the tangent space 
to $S_{t}$ at some point $p$ is the orthogonal of $p-\pi(p)$. 
Hence, $n(p) =  \frac{p-\pi(p)}{\sqrt{-\langle p-\pi(p) \mid p-n(p) \rangle}}$
is the normal vector to $S_{t}$ at $p$ pointing in the future. For this
reason, we call $n: S_{t} \rightarrow {\mathbb H}^{n-1}$ the \emph{Gauss map}.
Moreover, corollary $4.5$ expresses much more: 
$S(p+\epsilon n(p)) = S(p) + \epsilon$. Thus, from one
level set, let's say, $S_{1}$, and the Gauss map $n$ on it, 
we can reconstruct all other level sets.

\begin{prop}
\label{lip}
The Gauss map $n: S_{t} \rightarrow {\mathbb H}^{n-1}$, where
${\mathbb H}^{n-1}$ is equipped with its usual hyperbolic metric 
is $\frac{1}{t}$-Lipschitz.
\end{prop}

\preua
Select any $0 < \epsilon < t$. When $p$ describe $S_{t}$,
then $q(p) = p-\epsilon(p)n(p)$ describe $S_{t-\epsilon}$, and $n(q(p))=n(p)$. 
Since $T$ is $C^{1}$ and has spacelike fibers,
$S_{t}$ admits an induced riemannian metric, and for any $p$, $p'$
in $S_{t}$, the distance between them is the infimum of
$\sum_{i=0}^{N-1} \mid p_{i+1}-p_{i} \mid$ when 
$p_{0}, \ldots, p_{N}$ describe all the finite sequences in $S_{t}$
with $p_{0} = p$, $p_{N}=p'$. 
For such a sequence, we have:

\[ \mid p_{i+1} - p_{i} \mid^{2}  =  \mid q(p_{i+1}) - q(p_{i}) \mid^{2} + 
\epsilon^{2} \mid n(q(p_{i+1})) - n(q(p_{i})) \mid^{2} + 
2\epsilon\langle q(p_{i+1})-q(p_{i}) \mid n(q(p_{i+1})) - n(q(p_{i})) \rangle \]

The convexity of $S_{t-\epsilon}$ implies the positivity of
$\langle q(p_{i+1})-q(p_{i}) \mid n(q(p_{i+1})) - n(q(p_{i})) \rangle$.
Hence:

\[
\mid n(p_{i+1}) - n(p_{i}) \mid^{2}  =   \mid n(q(p_{i+1})) - n(q(p_{i})) \mid^{2}  
\leq  \frac{1}{\epsilon^{2}}  \mid p_{i+1} - p_{i} \mid^{2} \]

Since $\epsilon$ can be selected arbitrarly near $t$, the
proposition follows. \fin

\begin{cora}
Every $S_{t}$ is complete.
\end{cora}

\preua
Since $S_{t}$ is riemannian, completeness notions are all equivalent.
Let $s \mapsto p(s)$ be an incomplete geodesic in $S_{t}$
defined on $[0, s_{\infty}[$, parametrized by unit length.
According to proposition
\ref{lip}, $s \mapsto n(c(s))$ must converge
in ${\mathbb H}^{n-1}$ to some limit point
$n_{\infty}$. We select a coordinate system on \mink so that
this limit vector is $(1, 0, \ldots, 0)$ and such that the
Minkowski norm is $-dx_{0}^{2} + dx_{1}^{2} + \ldots + dx_{n-1}^{2}$.
For some small $\alpha$, every tangent vector $\partial_{t}c(t)$ with
$s-\alpha < s <s_{\infty}$ has $x_{0}$-component less than, 
let's say, $\frac{1}{2}$.
It follows that the orthogonal projection on $P_{0}=\{ x_{0} = 0 \}$ of
the geodesic $c$ has finite length for the usual euclidean metric of
$P_{0}$. Hence, $c(s)$ has a limit point $c_{\infty}$ in $P_{0}$.
We claim that the vertical line above $c_{\infty}$, like any timelike
line in \mink, intersects $\partial\Omega^{+}(\Lambda)$: 
indeed, it must enter in the future cone of every point in
$\Omega^{+}(\Lambda)$. Since  $\Omega^{+}(\Lambda)$ is
geodesically complete in the future, the vertical line
must thus intersect  $\Omega^{+}(\Lambda)$. On the other
hand, it cannot be entirely contained in $\Omega^{+}(\Lambda)$
since it intersects every degenerate hyperplanes, and thus,
every element of $\Lambda$.
 
From the geodesic completeness in the future of $\Omega^{+}(\Lambda)$,
it follows now that the vertical line intersects $S_{t}$ at one and only
one point $p_{\infty}$. The initial geodesic can then be completed
on $[0, s_{\infty}]$ by $c(s_{\infty}) = p_{\infty}$.\fin

Finally, as it is proved in \cite{bonsante}, lemma $4.9$, every $S_{t}$
is a Cauchy hypersurface for $\Omega(\Lambda)$. Thus:

\begin{prop} 
Future complete regular convex domains are globally hyperbolic
Cauchy-complete, admitting as Cauchy hypersurfaces the level
sets of the cosmological time function.\fin
\end{prop}

\rquea
We didn't systematically state all the similar results for 
past complete regular convex domains, but of course
they are true. 

We also could have written this section in a slightly different
way, by defining future (resp. past) regular convex domains as defined
by closed subsets of ${\mathcal J}^{+}$ (resp. ${\mathcal J}^{-}$)
(see remark \ref{ambigu}),
but it is useful for the next section to stress out that these closed subsets
arise from closed subsets in the same space $\mathcal J$.
\erquea

\subsection{Groups preserving regular convex domains and
Cauchy-hyperbolic spacetimes}

In the first part of this section, $\Gamma$ is a discrete subgroup of 
$\mbox{Isom}(M^{n})$. 

\begin{prop}
\label{agissons}
If $\Gamma$ preserves a future complete regular convex domain $\Omega$, 
then, the action of $\rho(\Gamma)$ on $\Omega$ is properly discontinuous.
If moreover $\Gamma$ is torsionfree, then this action is
free, and the quotient $\Gamma\backslash\Omega$ is
a Cauchy-complete semicomplete GH spacetime.
\end{prop}

\preua
According to the preceding \S, $\Omega$ is globally hyperbolic,
with a regular cosmological time function 
$T: \Omega \rightarrow ]0, +\infty[$. 
This function is $\Gamma$-invariant. The action of $\Gamma$ on every
level set $S_{t}$ is isometric for the induced riemanian
metric. It follows that the action on
$\Omega$ is properly discontinuous. 
Observe that an element of $\Gamma$ is trivial as soon
as its action on a level set $S_{t}$ is trivial.
Indeed, for any $x$ in such a $S_{t}$, $\gamma$
is an isometry of \mink which is trivial on the spacelike
hyperplane $T_{x}S_{t}$: since it preserves chronological
orientation, it follows that $\gamma$ is trivial.

Assume now that $\gamma$ admits a fixed point $x$ in $\Omega$.
Then, since its action on $S_{t}$ with 
$t = T(x)$ is isometric, and since $\Gamma$ is discrete,
it follows that $\gamma$ has finite order when restricted
to $S_{t}$. According to the above, $\gamma$ itself has
then finite order. Hence, if $\Gamma$ is torsionfree,
the action is free. The quotient space $\Gamma\backslash\Omega$
is then a spacetime, admitting the cosmological time function
induced by $T$ expressing it as a GH spacetime. The
level sets of this time function are Cauchy hypersurfaces
which are quotients $\Gamma\backslash S_{t}$:
the induced riemannian metric is complete.\fin

From now, $\Gamma$ denotes a discrete 
subgroup of $SO_{0}(1,n-1)$ without global fixed
point on $\overline{\mathbb H}^{n-1}$. 

Let $\rho: \Gamma \rightarrow {\bf M}^{n}$ such that
$L \circ \rho$ is the identity map. We say that $\rho$
is \emph{future admissible} (resp. past admissible) if 
$\rho(\Gamma)$ preserves a future complete (resp. past complete)
regular convex domain.

According to Proposition \ref{agissons},
admissible representations produce Cauchy-complete
GH spacetimes. We want to associate to every such representation
a \emph{maximal} GH spacetime; it will be provided by the
following proposition.

\begin{prop}
\label{admlox}
For any admissible representation $\rho$, the closure of
the set of repulsive fixed
points in ${\mathcal J}$ of loxodromic elements of $\rho(\Gamma)$
is a $\rho(\Gamma)$-invariant subset $\Lambda(\rho)$ 
contained in any closed $\rho(\Gamma)$-invariant subset
of ${\mathcal J}$. 
\end{prop}

\preua
Consider an admissible representation $\rho$, preserving a 
closed subset $\Lambda$. The fibration 
$\delta: {\mathcal J} \rightarrow {\mathcal S}$ is 
$\Gamma$-equivariant, where the $\Gamma$-action on 
$\mathcal S$ is its usual conformal action on the sphere.
The dynamic of (nonelementary) discrete subgroups of 
$SO_{0}(1,n-1)$ on the conformal sphere
is well-known: there is a closed subset $\widehat{\Lambda}$,
the limit set, contained in every $\Gamma$-invariant
closed subset - in particular, in the closure of
$\delta(\Lambda)$.
It is also well-known that $\widehat{\Lambda}$ is the closure
of repulsive fixed points of loxodromic elements of
$\Gamma$. Let $\gamma$ such a loxodromic element
of $\Gamma$, and $v_{0}$ the element of $\mathcal S$
corresponding to the repulsive fixed point of $\gamma$.
We consider here $v_{0}$ as a future oriented lightlike
vector in \mink with $N$-norm $1$, then, $v_{0}$ is
a $\gamma$ eigenvector, with eigenvalue $0 < \lambda < 1$. 
The bassin of repulsion of $x_{0}$ for $\gamma$ contains some
open subset of the form $\delta^{-1}(U)$, where
$U$ is some open neighborhood of $v_{0}$ in
$\mathcal S$. Then, $U$ contains some element
of $\delta(\Lambda)$, hence, $\delta^{-1}(U)$
contains some element $x$ of $\Lambda$.
Then, negative iterates of $x$ under
$\gamma$ converge towards $x_{0}$.

It follows as required that every repulsive fixed point of loxodromic
elements of $\rho(\Gamma)$ belongs to $\Lambda$. \fin

\begin{defina}
\label{definitioncauchyhyp}
A Cauchy-hyperbolic spacetime is the quotient of 
$\Omega^{\pm}(\Lambda(\rho))$
by $\rho(\Gamma)$ for any admissible representation 
$\rho: \Gamma \rightarrow \mbox{Isom}({\bf M}^{n})$.
\end{defina}

\rquea
It follows from proposition \ref{admlox} and \S \ref{maxi} below
that any Cauchy-complete GH spacetime with holonomy group $\rho(\Gamma)$
can be isometrically embedded in the associated Cauchy-hyperbolic
spacetime.
\erquea

\subsection{Admissible and non-admissible representations}
Let $\Gamma$ be a nonelementary discrete subgroup of
$SO_{0}(1,n-1)$. We consider 
the space of representations $\rho: \Gamma \rightarrow \mbox{Isom}({\bf
M}^{n})$ admitting the identity map as linear part.
For such a $\rho$, the translation part $\tau(\gamma)$
in the expression 
$\rho(\gamma)x = \gamma x + \tau(\gamma)$ defines a map 
$\tau: \Gamma \rightarrow {\mathbb R}^{1,n-1}$ which satisfies
$\tau(\gamma\gamma') = \tau(\gamma) + \gamma\tau(\gamma')$,
i.e. which is a $1$-cocycle for the $\Gamma$-module
${\mathbb R}^{1,n-1}$. We denote by $Z^{1}(\Gamma, {\mathbb R}^{1,n-1})$
the space of cocycles.

Two such cocycles defines representations
conjugate by a translation of \mink if and only if they
differ by a coboundary, i.e., a cocycle of the form
$\tau(\gamma) = \gamma v - v$ for some fixed $v$. We denote
by $B^{1}(\Gamma, {\mathbb R}^{1,n-1})$ the space of coboundaries.

In other words, the space of representations $\rho$
as above up to conjugacy by translations 
is parametrized by the quotient of
$Z^{1}(\Gamma, {\mathbb R}^{1,n-1})$ by $B^{1}(\Gamma, {\mathbb R}^{1,n-1})$,
i.e., the twisted cohomology
space $H^{1}(\Gamma, {\mathbb R}^{1,n-1})$, the $\Gamma$-module
structure on ${\mathbb R}^{1,n-1}$ being given by the fixed linear part 
$\Gamma$. Observe that $H^{1}(\Gamma, {\mathbb R}^{1,n-1})$
is naturally equipped with a structure of linear space.

Here, we want to describe the set of \emph{future admissible\/}
representations, i.e. the domain ${\mathcal T}^{+}$
in $H^{1}(\Gamma, {\mathbb R}^{1,n-1})$
corresponding to representations preserving a future complete
regular convex domain.
According to proposition \ref{admlox}, this is exactly the set of cocycles
corresponding to representations for which the future complete convex
set defined by repulsive fixed points of loxodromic elements is not
empty.

Similarly, we can define the domain ${\mathcal T}^{-}$ 
correponding to representations preserving some past complete
regular domains, but observe that the conjugacy 
by $-id$ induces a transformation on 
$H^{1}(\Gamma, {\mathbb R}^{1,n-1})$, which is nothing but the antipodal map,
and exchanges ${\mathcal T}^{+}$, ${\mathcal T}^{-}$.
In other words, we have 
${\mathcal T}^{-} = -{\mathcal T}^{+}$, 
hence we restrict our study to future complete regular domains.

The conjugacy by some positive homothety $\lambda id, \;\;\; \lambda>0$
preserves the linear part too and induces on 
$H^{1}(\Gamma, {\mathbb R}^{1,n-1})$ a positive homothety.
It clearly preserves ${\mathcal T}^{+}$ which is thus a cone.

\begin{lema}
${\mathcal T}^{+}$ is a convex cone.
\end{lema}

\preua
Let $[\tau_{1}]$, $[\tau_{2}]$ be two elements of
${\mathcal T}^{+}$. Let's denote $\Lambda_{i}$ the closure
in ${\mathcal J}$ of the set of repulsive fixed points
of loxodromic elements of $\Gamma$ for the representation
$\rho_{i}$ associated to $\tau_{i}$. Then,
$\Omega^{+}(\Lambda_{i})$ for $i=1,2$ are nonempty
future complete regular domains, and thus, have nonempty
intersection. Let $p$ be an element of 
$\Omega^{+}(\Lambda_{1}) \cap \Omega^{+}(\Lambda_{2})$.

For any loxodromic element $\gamma$ of $\Gamma$,
we denote by $(v(\gamma), s_{i}(\gamma))$ its unique
repulsive fixed point in 
${\mathcal J} \approx {\mathcal S} \times \mathbb R$ for
the representation $\rho_{i}$. Observe that $v(\gamma)$ is indeed
the same for $i=1,2$, since it is the future oriented $\gamma$-eigenvector
with $N$-norm $1$ associated to the eigenvalue $\lambda$ with
absolute value less than $1$. We have (see \S \ref{penrose}):

\begin{eqnarray*} 
( v(\gamma) , s_{i}(\gamma) ) & = & 
\rho_{i}(\gamma)( v(\gamma) , s_{i}(\gamma) )\\
 & = & ( v(\gamma) , \lambda^{-1}s_{i}(\gamma)+\lambda^{-1}
\langle \tau_{i} \mid \lambda v(\gamma) \rangle)
\end{eqnarray*}

Hence: 
\[ s_{i}(\gamma) = -\frac{\langle \tau_{i} \mid v(\gamma)
\rangle}{\lambda^{-1} - 1} \]

Now, since $p$ belongs to 
$\Omega^{+}(\Lambda_{1}) \cap \Omega^{+}(\Lambda_{2})$, for any
loxodromic $\gamma$ and any $i$, we have 
$\langle v(\gamma) \mid p \rangle - s_{i}(\gamma) < 0$, 
i.e.:

\[ (\lambda^{-1}-1)\langle v(\gamma) \mid p \rangle + \langle \tau_{i} \mid
v(\gamma) \rangle  < 0 \]

Clearly, this last expression is still valid
if the term $\tau_{i}$ is replaced by any $\alpha\tau_{1}+(1-\alpha)\tau_{2}$
with $0 \leq \alpha \leq 1$. It follows that representations associated
to $[\alpha\tau_{1}+(1-\alpha)\tau_{2}]$ admits as set of
repulsive fixed points
in ${\mathcal J}$ a subset for which the associated future complete convex 
contains $p$, i.e., is not empty. This means precisely that these
representations all belong to ${\mathcal T}^{+}$. \fin

\rquea
The intersection ${\mathcal T}^{+} \cap {\mathcal T}^{-}$
is a linear subspace since it is a convex cone stable by antipody.
Observe that this intersection correspond to representations 
preserving a \emph{compact} subset of ${\mathcal J}$
(cf. corollary \ref{futetpast}).
\erquea

\subsubsection{Non admissible cocycles}

\begin{prop}
\label{tri}
There his a (nonuniform) lattice $\Gamma$ of $SO_{0}(1,2)$ for which
${\mathcal T}^{+}$ is reduced to $\{ 0 \}$.
\end{prop}

\preua
Consider the $3$-punctured sphere, i.e. the riemann surface
of genus $0$ with $3$ cusps (there is only one such 
riemann surface). It is the quotient of the Poincar\'e disk
by a fuchsian group $\Gamma$, which is isomorphic to
the free group of rank $2$, generated by three parabolic
elements $a$, $b$ and $c$ of $SO_{0}(1,2)$ satisfying
the relation $abc = id$. There is another description
of $\Gamma$: consider an ideal triangle in the Poincar\'e 
disk, and the group generated by reflexions around edges
of this triangle. This group contains an index $2$ 
subgroup, the orientation preserving elements, which
is nothing but the group $\Gamma$. It is clear 
from this last description that there is an elliptic
element $R$ of $SO_{0}(1,2)$ of order $3$ - the rotation 
permuting the ideal vertices of the initial ideal
triangle - such that the conjugacy by $R$ cyclically
permutes $a$, $b$ and $c$.
Denote by $\alpha$, $\beta$ and $\kappa$ the
unique isotropic vectors fixed respectively
by $a$, $b$ and $c$ (they are the vertices of the initial
ideal triangle). Of course, $\kappa = R(\beta) =
R^{2}(\alpha)$.

Consider now $H^{1}(\Gamma, {\mathbb R}^{1,2})$.
For any cocycle $\tau$, every $\tau(\gamma)$ 
can be computed from $\tau(a)$ and $\tau(b)$, since
$a$ and $b$ generate $\Gamma$. Hence, cocycles form
a $6$-dimensional linear space naturally identified
with ${\mathbb R}^{1,2} \times {\mathbb R}^{1,2}$.
Coboundaries are the image of the map from
${\mathbb R}^{1,2}$ into ${\mathbb R}^{1,2} \times {\mathbb R}^{1,2}$
which associates to $x$ the pair $(ax-x, bx-x)$.
This map is injective, since $a$, $b$ have no (nontrivial) common fixed points,
hence, the image is $3$-dimensional: 
$H^{1}(\Gamma, {\mathbb R}^{1,2})$ has dimension $3$.

It will be useful later to represent elements of
$Z^{1}(\Gamma, {\mathbb R}^{1,2})$ by triples
$(\tau(a), \tau(b), \tau(c))$ satisfying the
cocycle relation $\tau(a) + a\tau(b) + ab \tau(c) = 0$.

Now, if $\tau$
belongs to ${\mathcal T}^{+}$, the Minkowski
isometries associated by this cocycle to $a$, $b$
and $c$ all preserve a regular convex
domain, and thus, a spacelike complete
hypersurface. According to \S \ref{ach} below,
these isometries are not transverse, i.e.,
the translations vectors $\tau(a)$, $\tau(b)$
and $\tau(c)$ are respectively orthogonal to
$\alpha$, $\beta$ and $\kappa$. Denote by $E$
the space of cocycles satisfying these
orthogonality conditions: we have just proved that
$E$ contains ${\mathcal T}^{+}$. Observe
that $E$ is defined by linear conditions on
$H^{1}(\Gamma, {\mathbb R}^{1,2})$; it is 
therefore a linear subspace. Moreover, since
$\tau(a)$, $\tau(b)$ can be selected non orthogonal
to $\alpha$, $\beta$, the codimension of $E$ is at least $2$.
At first glance, one could think that the third
condition ``$\tau(c)$ is orthogonal to $\kappa$''
immediatly implies that $E$ has codimension $3$, i.e.,
is trivial. But this is not clear since it
could be true that by some extraordinary miracle,
$\tau(c)$ is orthogonal to $\kappa$ as soon as
$\tau(a)$, $\tau(b)$ are themselves orthogonal
to $\alpha$, $\beta$. In this case, $E$ would
be a line.

To prove that this miracle does not occur, we consider the rotation
$R$ above: the conjugacy by $R$ induces naturally
an action on $Z^{1}(\Gamma, {\mathbb R}^{1,2})$:
if $\tau$ is a cocycle, 
$R(\tau)(\gamma) = R\tau(R^{-1}\gamma R)$. 
When we express elements of $Z^{1}(\Gamma, {\mathbb R}^{1,2})$
by triples $(\tau(a), \tau(b), \tau(c))$, the iteration
under $R$ maps such an element on 
$(R\tau(c), R\tau(a), R\tau(b))$. 
Now, we observe that this action has order $3$ 
(exactly order $3$ since
it is not trivial). Hence, as any nontrivial automorphism
of $3$-dimensional spaces of order $3$, the action induced
by $R$ one the quotient $H^{1}(\Gamma, {\mathbb R}^{1,2})$
has a line of fixed points $\Delta$, and a plane on which it acts
as a rotation with angle 
$\frac{2\pi}{3}$. In particular,
$\Delta$ is the unique line globally preserved by $R$.
Hence, if the miracle imagined above occurs, $E$, which is 
obviously $R$-invariant, must be equal to $\Delta$.

The next step is the identification of $\Delta$:
if a cocycle $\tau$ represents a fixed point of $R$ in
$H^{1}(\Gamma, {\mathbb R}^{1,2})$, it is cohomologous to
$\frac{\tau + R(\tau) + R^{2}(\tau)}{3}$. Hence,
elements of $H^{1}(\Gamma, {\mathbb R}^{1,2})$ fixed by $R$
are all represented by triples
$(\tau(a), \tau(b), \tau(c))$ satisfying
$\tau(a) = R\tau(c)$, $\tau(b)=R\tau(a)$, 
$\tau(c)=R\tau(b)$. The cocycle property of such
a triple reduces to:

\[
\tau(a) + aR\tau(a) + abR^{2}\tau(a)   =  0 \]

i.e.:
\[
\tau(a) + aR\tau(a) + (aR)^{2}\tau(a)  =  0
\]

Thus, the element $aR$ of $SO_{0}(1,2)$ is of special
interest. Observe that it has order $3$. Indeed:

\begin{eqnarray*}
(aR)^{3} & = & a(Ra)(Ra)R \\
         & = & ab(Rb)R^{2} \\
	 & = & abcR^{3} \\
	 & = & id
\end{eqnarray*}

It is also nontrivial, since $aR(\kappa) = \alpha$.
Hence, in ${\bf M}^{3}$, it admits a timelike line of fixed points,
and the orthogonal to this line of fixed points
is a spacelike plane $P_{0}$. Now, we observe that the cocycle
property stated above for $(\tau(a), R\tau(a), R^{2}\tau(a))$
means precisely that $\tau(a)$ must belong to $P_{0}$!

Hence, the line of fixed points $\Delta$ is a quotient
of $P_{0}$ by elements $\tau(a)$ representing $R$-invariants
\emph{coboundaries}. These coboundaries form a subspace of
$P_{0}$ of dimension $1$, since $\Delta$ has dimension $1$,
and $P_{0}$ has dimension $2$. We claim that these
coboundaries are precisely $\tau(a) = ax-x$ where $x$ is a
fixed point of $R$. Indeed:

- the triple $(ax-x, bx-x, cx-x)$ is $R$-invariant:
$R(ax-x) = Rax-Rx = bRx - x = bx-x$, and similarly,
$R(bx-x) = cx-x$,

- the vector $ax-x$ belongs to $P_{0}$. Indeed:
\begin{eqnarray*}
(ax-x) + aR(ax-x) + (aR)^{2}(ax-x) & = & ax-x + aRax-ax + aRaRax - aRax\\ 
                                   & = & (aR)^{3}x - x \\
				   & = & 0
\end{eqnarray*}

Now, we observe that the $ax-x$ are exactly the vectors in $P_{0}$
orthogonal to $\alpha$. In other words, the $R$-fixed triples 
$(\tau(a), R\tau(a), R^{2}\tau(a))$ representing elements in $E$
are all coboundaries. It means that the intersection between
$E$ and $\Delta$ is reduced to $\{ 0 \}$.
In particular, $E \neq \Delta$: as observed previously,
it implies that $E$ is reduced to $\{ 0 \}$.\fin

\rquea
The proof above highly relies on the specific symmetric
properties of the group. The identification of ${\mathcal T}^{+}$
for geometrically finite Kleinian groups remains an interesting challenge,
even in dimension $2+1$.
\erquea

\subsubsection{Convex cocompact Kleinian groups}
\label{suli}
 
We first recall some well-known facts on Kleinian group. See e.g.
\cite{ratcliffe}.
Let $\Gamma$ be a such a Kleinian group, i.e. a finitely generated
(nonelementary) discrete subgroup of $SO_{0}(1,n-1)$.
Let $\widehat{\Lambda}$ be its limit set in the conformal sphere
${\mathcal S} \approx \partial{\mathbb H}^{n-1}$, and let
$C(\widehat{\Lambda})$ be the convex hull of $\widehat{\Lambda}$ in
${\mathbb H}^{n-1}$. This is the minimal $\Gamma$-invariant
closed convex subset of ${\mathbb H}^{n-1}$.

\begin{defina} 
$\Gamma$ is \emph{geometrically finite} if, for any $\epsilon>0$,
the quotient of the $\epsilon$-neighborhood of $C(\widehat{\Lambda})$ by
$\Gamma$ has finite volume. If moreover this quotient is
compact, then $\Gamma$ is \emph{convex cocompact.}
\end{defina}

The main aspect of convex cocompact groups we will use
is the \emph{hyperbolic} character of their dynamics around
the limit sets. Namely (see \cite{sullivan}, \S $9$), there is
a finite symetric generating set $G$ for $\Gamma$ of
loxodromic elements $\gamma_{1}, \ldots , \gamma_{k}$
such that:

(i) for some fixed round metric on ${\mathcal S}$ in the natural conformal class,
if $U_{i}$ denotes the (open) domain of the sphere where $g_{i}$ is expanding,
then the union of the $U_{i}$ covers $\widehat{\Lambda}$,

(ii) there is an uniform $N$ such that,
for any $x_{0}$ in $\widehat{\Lambda}$, and for any pair of sequences
$\gamma_{i_{0}}, \gamma_{i_{1}}, \ldots $ and $\gamma_{j_{0}},
\gamma_{j_{1}}, \ldots$ satisfying 
$x_{n+1} = \gamma_{i_{n}}x_{n} \in U_{i_{n+1}}$ and 
$x'_{n+1} = \gamma_{j_{n}}x'_{n} \in U_{j_{n+1}}$, then 
$(\gamma_{j_{n}} \ldots \gamma_{j_{0}})(\gamma_{i_{n}} \ldots
\gamma_{i_{0}})^{-1}$ is equal to a product 
$\gamma_{k_{l}} \ldots \gamma_{k_{0}}$ with $l \leq N$.

Actually, D. Sullivan in \cite{sullivan} states the hyperbolic property
of convex cocompact Kleinian groups only in the $3$-dimensional case,
but his proof applies for the general case in any dimension:
indeed, the proof in \cite{sullivan} relies on the fact
that the limit set admits only \emph{conical} limit points,
and this fact remains true in any dimension (in fact, it is another
equivalent definition of convex cocompact groups, see \cite{ratcliffe})).
Anyway, any reader acquainted with dynamical systems theory will
recognize this hyperbolic property on $\widehat{\Lambda}$ as the
hyperbolic property of the geodesic flow on $\Gamma\backslash{\mathbb H}^{n}$
around the compact invariant set formed by the geodesics lying entirely in 
$\Gamma\backslash C(\widehat{\Lambda})$.

Observe that this hyperbolic
property extends directly to the action of 
$\Gamma \subset SO_{0}(1,n-1) \subset \mbox{Isom}({\bf M}^{n})$
on ${\mathcal J}$, where the compact invariant subset
is $\widehat{\Lambda} \times \{ 0 \}$ (here, the action
on \mink under consideration is the action associated to the trivial cocycle).
Indeed, the property $(i)$ follows from proposition \ref{laoucapulse}
and the expanding property of $\Gamma$ on $\mathcal S$, and the property
$(ii)$ follows directly from its version on $\mathcal S$.

Hence, and as observed in the last remark page $259$ of \cite{sullivan},
Theorem $II$ of \cite{sullivan} can be applied: the action
of $\Gamma$ on $\widehat{\Lambda} \times \{ 0 \}$ is
structurally stable. In other words, for any action of $\Gamma$
on ${\mathcal J}$ $C^{0}$-near the initial hyperbolic action,
there is a compact invariant subset, on which the action restricts
as an action topologically conjugate to the restriction of
$\Gamma$ on $\widehat{\Lambda} \times \{ 0 \}$. 
This is true in particular for actions associated to nontrivial cocycles
in $Z^{1}(\Gamma, {\mathbb R}^{1,n-1})$ which are sufficiently small,
i.e., for which the $\tau(\gamma_{i}) \;\;\; (1 \leq i \leq k)$ have
small $N$-norm. These small cocycles correspond to actions on
${\mathcal J}$ preserving a compact subset; therefore, they
belong to ${\mathcal T}^{+} \cap -{\mathcal T}^{+}$. Since
${\mathcal T}^{+}$ is a cone:

\begin{thma}
For convex cocompact Kleinian groups, the convex cone 
${\mathcal T}^{+}$ is the whole $H^{1}(\Gamma, {\mathbb R}^{1,n-1})$.
\fin
\end{thma}

\section{Around Bieberbach's Theorem}

\subsection{Non discrete linear parts and Auslander's Theorem}
\label{bieb}
 
We will use the following theorem, refinement by Y. Carri\`ere and
F. Dal'bo of a theorem by L. Auslander (see \cite{cardal}, Theorem $1.2.1$):

\begin{thma}
\label{ausl}
Let $\Gamma$ be a discrete subgroup of $\mbox{Aff}(n, {\mathbb R})$.
Then, the identity component $G_{0}$ of
the closure of $L(\Gamma)$ in $GL(n, {\mathbb R})$ is nilpotent.
In particular, $\Gamma_{nd} = L^{-1}(G_{0} \cap L(\Gamma))$
is nilpotent.
\end{thma}

This theorem fits perfectly with (see \cite{cardal},
Proposition $1.2.2$):

\begin{thma}
\label{unip}
Let $\Gamma$ be a nilpotent group of affine transformations of
${\mathbb R}^{n}$. Then, there exist a maximal $\Gamma$-invariant 
affine subspace ${\bf U}$
of ${\mathbb R}^{n}$ such that the restriction of the action
of $\Gamma$ to ${\bf U}$ is unipotent. Moreover, this maximal
unipotent affine subspace is unique.
\end{thma}

As an application of these theorems, it is easy for example to
recover the following version of Bieberbach's Theorem:

\emph{If $\Gamma$ is discrete group of isometries of the euclidean
space of dimension $n$, it contains a finite index free abelian subgroup of finite
rank.\/} 

The proof goes as follows: if $L(\Gamma) \subset SO(n)$ is discrete, it is finite;
the kernel of $L_{\mid \Gamma}$ is a finite index subgroup
of translations. If $L(\Gamma)$ is not discrete, we consider the
unique maximal unipotent affine subspace associated to $\Gamma_{nd}$:
$\Gamma_{nd}$ is a finite index subgroup of $\Gamma$, and the
kernel of the restriction morphism $\Gamma \rightarrow \mbox{Isom}({\bf U})$
is a finite index subgroup of $\Gamma_{nd}$. The proof is completed by
the observation that unipotent elements of $\mbox{Isom}({\bf U})$
are translations.

The most famous version of Bieberbach's theorem is maybe the
following corollary:

\emph{If $\Gamma$ is a cristallographic group, i.e. a discrete
group of isometries such that ${\mathbb R}^{n}\slash\Gamma$ is compact,
then it contains a finite index subgroup of translations.\/}

This complement arises from a homological argument:
since ${\mathbb R}^{n}\slash\Gamma$ is compact, any finite index
subgroup has a non-trivial Betti-number $b_{n}$, but since 
there is a finite index subgroup of $\Gamma$ acting as translations
on ${\bf U}$, this finite index subgroup has trivial Betti numbers $b_{i}$
for $i \geq \mbox{dim}(\bf U)$.

\subsection{Twisted product by euclidean manifolds}
\label{produi}

Truely speaking, it is easy to provide a more precise version of
Bieberbach's theorem in the non compact case, but requiring a
new notion: let $M$ be any manifold, and denote by $\Gamma$ its
fundamental group. Let $N$ be any flat euclidean manifold,
and $r: \Gamma \rightarrow \mbox{Isom}(N)$ be any morphism.
It is well-known how to produce from such a data a 
locally trivial bundle $q: B \rightarrow M$ with fibers
diffeomorphic to $N$: if $\widetilde{M}$ is the universal covering
of $M$, consider the diagonal action 
$(\tilde{x}, y) \mapsto (\gamma\tilde{x}, r(\gamma)(y))$ of $\Gamma$
on the product $\widetilde{M} \times N$. Since the action of
$\Gamma$ on the first component $\widetilde{M}$ is free and properly
discontinuous, the same is true for its action on $\widetilde{M} \times N$.
Denote by $B$ the quotient space, and $q: B \rightarrow M$ the
map induced by the first projection map. Clearly, $q$ is a locally
trivial fibration, with fiber $N$ as claimed above. It is called
the \emph{suspension over $M$ with monodromy $r$.\/}

Now, if $M$ is a flat euclidean manifold, the product metric on 
$\widetilde{M} \times N$ is locally euclidean 
(since $N$ is assumed here locally euclidean) and is preserved
by the $\Gamma$-action defined above. Therefore, the quotient
space $B$ inherits itself a locally euclidean structure, canonically
defined from the initial locally euclidean
structures on $M$ and $N$. $B$, equipped with this euclidean structure,
is called the \emph{twisted product of $M$ by $N$ with monodromy
$r$.\/}

The argument in the preceding section actually show that
\emph{any locally euclidean manifold is a twisted product of
a flat torus $\Gamma \backslash{\bf U}$ by a euclidean
linear space ${\mathbb R}^{n}/{\bf U}$ with holonomy
in $SO({\mathbb R}^{n}/{\bf U})$.\/}

It should be clear to the reader that a similar procedure
can be defined when $M$ is a flat lorentzian manifold;
the result is then a canonical flat lorentzian structure
on $B$, which is called once more \emph{the twisted product
of $M$ by the euclidean manifold $N$, with monodromy
$r$.\/} When $N$ is an euclidean linear space and $r$
takes value in linear isometries (i.e., $r(\Gamma)$ admits
a global fixed point in $N$), then we say that
$B$ is a \emph{linear twisted product over $M$.}

\rquea
If $t: M \rightarrow {\mathbb R}$ is a proper time function 
expressing $M$ as a GH spacetime, the composition $t \circ q$
satisfies the same properties. Hence, \emph{the twisted products
by euclidean manifolds of flat GH spacetimes are still
flat GH spacetimes.\/} The inverse statement, namely the fact 
that $M$ is GH as soon as 
$B$ is, even if far from obvious now, follows from theorem \ref{caucomp}. 

In the same vein, twisted products over Cauchy-complete
GH spacetimes are still Cauchy-complete, but, of course,
such a procedure preserves Cauchy-compactness if and only if
the euclidean manifold $N$ is closed (i.e. finitely covered
by a flat torus).
\erquea

\rquea
Since here we don't worry about finite index phenomena, the Bieberbach's
theorem as stated above gives a completely satisfactory description
of euclidean manifolds, and their isometry groups $\mbox{Isom}(N)$
are easily described. For example, when $N$ is a flat torus,
up to finite index, the representation 
$r: \Gamma \rightarrow \mbox{Isom}(N)$ can be assumed as taking
value in translations on $N$.

We should also point out that all twisted products we consider
in this work are either linear twisted products, or twisted
products by flat torii.
\erquea

The following lemma will be useful to recognize twisted products,
and its proof should be obvious to the reader:

\begin{lema}
\label{suspendre}
Let $B$ be a flat spacetime, 
quotient of an open subset $\Omega$ of \mink by a 
group of isometries acting freely and properly discontinuously.
Assume that $\Gamma$ preserves a timelike affine subspace $\bf U$,
and that $\Omega$ is invariant by translations by vectors contained
in the orthogonal linear space ${\bf U}^{\perp}$. Then, the
quotient of ${\bf U} \cap \Omega$ by $\Gamma$ is a 
flat spacetime $M$, and $B$ is the linear twisted product of
$M$ by the euclidean space ${\bf U}^{\perp}$.\fin
\end{lema}

\section{Classification of isometries}
\label{class}

Let $g$ be a isometry of ${\bf M}^{n}$. We denote by $L(g)$
its linear part, $p$ the element $L(g)-id$ of $SO_{0}(1,n-1)$,
and $\tau$ the translation part: $g(x)=L(g)(x)+\tau$.

Observe that:

\[ (\ast) \;\;\;\;\; \langle x \mid p(y) \rangle + \langle p(x) \mid y \rangle +
\langle p(x) \mid p(y) \rangle = 0 \]

In particular, if we denote by $I$ the image of $p$, the kernel
of $p$ is the orthogonal $I^{\perp}$.

Moreover, $\tau$ is defined modulo $I$, since it can be modified
through conjugacies by translations.

Up to conjugacy, every isometry is of the following form:

\begin{itemize}

\item \emph{Elliptic:\/} $L(g)$ preserves some euclidean norm on \mink.
Then, since ${\bf M}^{n} = I \oplus I^{\perp}$,
$\tau$ can be selected in $I^{\perp}$, i.e., fixed by $L(g)$.
Observe that in this case, $I$ is spacelike, except if $R$ is the 
identity map. Indeed, $L(g)$ considered as an isometry of
${\bf H}^{n-1}$ has a fixed point there, proving that
$I^{\perp}$ contains a timelike element. The claim follows.

\item \emph{Hyperbolic:\/} of the form:
\[ \left(\begin{array}{cccccc}
         ch(\zeta) & sh(\zeta) & 0 & 0 & \ldots & 0 \\
         sh(\zeta) & ch(\zeta) & 0 & 0 & \ldots & 0 \\
            0      &   0       & 1 & 0 & \ldots & 0 \\
            \ldots & \ldots & \ldots & \ldots & \ldots & \ldots \\
            0      &  0     &  0 & 0 & \ldots & 1
            \end{array}\right) \]
In this case, $\tau$ can be assumed in $I^{\perp}$. 
          
\item \emph{Unipotent:\/} $p$ is nilpotent. There is a trichotomy
on this case, depending on the position of $\tau$ with respect
to $I$ and $I^{\perp}$.         

\item \emph{Loxodromic:\/} of the form $R\circ A$, where $R$ is
elliptic, $A$ hyperbolic, and $R \circ A = A \circ R$.
Once more, $\tau$ should be assumed fixed by $L(g)$.

\item \emph{Parabolic:\/} of the form $R\circ A$, where $R$ is
elliptic, $A$ unipotent, and $R \circ A = A \circ R$.
                        
\end{itemize}

We have to understand better the parabolic case; in particular:

\subsection{Unipotent linear parts}

Observe that in this case, $J = I \cap I^{\perp}$ is a nontrivial
isotropic space, therefore, it is generated by a single isotropic element
$v_{0}$. 

Obviously, $p$ preserves $J^{\perp}$; it induces a nilpotent
transformation $\bar{p}$ on $J^{\perp}/J$. On the other hand,
the norm on ${\bf M}^{n}$ induces an euclidean norm on $J^{\perp}/J$
which is preserved by $id+\bar{p}$. It follows that
$\bar{p}$ is the zero map. Hence,
the image of $p^{2}$ is contained in $J$, since $I \subset J^{\perp}$.

Consider now an element $v_{1}$ of \mink such that 
$p(v_{1})=v_{0}$:

\[ 0=2\langle v_{1} \mid p(v_{1}) \rangle + \langle p(v_{1}) \mid 
p(v_{1}) \rangle =
2\langle v_{1} \mid v_{0} \rangle + \mid v_{0} \mid^{2} \]

It follows that $v_{1}$ must belong to $J^{\perp}$. 

Select now an element $v_{n}$ of \mink \emph{non\/}-orthogonal to $v_{0}$:
$v'_{1}=p(v_{n})$ is a non-trivial
element of $I \setminus J$. In particular, its norm is positive.

Eventually:

\[ 0=\langle p(v'_{1}) \mid v'_{1} \rangle + \langle v'_{1} \mid v'_{1} \rangle
+ \langle p(v'_{1}) \mid v_{n} \rangle \]

Since the middle term $\mid v'_{1} \mid^{2}$ is nonzero, $p(v'_{1})$
is not null.  Since the image of $p^{2}$ is contained in $J$,
$p(v'_{1}) = p^{2}(v_{n})$ generates $J$. 
Thus, we can eventually define $v_{0}$ as $p^{2}(v_{n})$,
and $v_{1}$ as $v'_{1}=p(v_{n})$. $L(g)$ is then represented , in
a non-orthogonal frame, by the matrix:

\[ \left(\begin{array}{cccccc}
   0 & 1 & 0 & \ldots & 0 & 0 \\
   0 & 0 & 0 & \ldots & 0 & 1 \\
   \ldots & \ldots & \ldots & \ldots & 0 & 0 \\
   0 & \ldots & \ldots & \ldots & \ldots & 0
   \end{array}\right)
   \]
   
For the translation vector $\tau$, we have three cases to consider:

- either it belongs to $I$, in which case we can assume it to be $0$;
we call this case the \emph{linear case.\/}

- either it belongs to $J^{\perp} \setminus I$, in which case
we can assume that it belongs to $I^{\perp}$, since 
$J^{\perp} = I + I^{\perp}$: we call it the \emph{tangent case.\/}

- either it does not belong to $J^{\perp}$: we call this
case the \emph{transverse case.\/}

\rquea
In the $3$-dimesional case, we have the equality $I = J^{\perp}$,
hence, the tangent case does not occur.
\erquea

\subsection{Parabolic linear parts}

We consider here the case $L(g) = R \circ (id + p) = (id+p) \circ R$
with $p$ nilpotent and $R$ a non-trivial elliptic element of \son.

The preceding section characterizes $p$. We observe that $I$, $J$
and their orthogonals
must be $R$-invariant. Since $R$ belongs to the orthochronous 
component, $R(v_{0}) = v_{0}$. Being elliptic, its restriction to 
$I = \langle v_{0}, v_{1} \rangle$ is the identity map.

Therefore, $R \circ p = p$. Denote by $\mathcal I$ the image of
$R-id$; fixed points of $R$ are the elements of ${\mathcal I}^{\perp}$.
The claim above implies the inclusions $I \subset {\mathcal I}^{\perp}$,
${\mathcal I} \subset I^{\perp}$. Observe also that
the orthogonal of $J^{\perp}$ by any $R$-invariant euclidean norm
is a $R$-invariant line: it is contained in 
the fixed point set of $R$, i.e., in ${\mathcal I}^{\perp}$.
In other words, $v_{n}$ can be selected $R$-invariant.

Now, we observe that $\tau$ is defined modulo $\mbox{Im}(R-id+p)$.
Any element $w$ of \mink is a sum $w=w_{1}+w_{2}$
with $w_{1}$ in $\mathcal I$ and $w_{2}$ in ${\mathcal I}^{\perp}$.
Then, $(R-id+p)(w_{1}) = (R-id)(w_{1})$ since $p(\mathcal I)=0$,
and $(R-id+p)(w_{2}) = p(w_{2})$. We deduce that 
$\mbox{Im}(R-id+p)$ contains $I$ and $\mathcal I$:
the first claim because ${\mathcal I}^{\perp}$ contains 
$I+\langle v_{n} \rangle$, and the second claim because 
$(R-id){\mathcal I}={\mathcal I}$. The reverse inclusions
being obvious ($p(w_{2}) \in I$ and $(R-id)(w_{1}) \in {\mathcal I}$),
we obtain that $\mbox{Im}(R-id+p) = I \oplus {\mathcal I}$,
the sum being direct and orthogonal.

Hence, $\tau$ can be assumed belonging to ${\mathcal I}^{\perp}$,
i.e., fixed by $R$. Inside ${\mathcal I}^{\perp}$, we can select $\tau$ 
modulo $I$: in particular, it can be assumed orthogonal to $v_{1}$. 

Once more, we distinguish three cases:

- the \emph{linear case\/} if $\tau$ belongs to $I$

- the \emph{tangent case\/} if $\tau$ belongs to $J^{\perp} \setminus I$

- the \emph{transverse case\/} if $\tau$ does not belong to $J^{\perp}$.

Observe that a parabolic element is nontransverse if and only if 
it preserves every affine hyperplane with direction $J^{\perp}$.

\section{Achronal domains of isometries}
\label{ach}
For any isometry $g$, we define its achronal domain $\Omega_{g}$ by:

\[ \Omega_{g} = \{ x \in {\bf M}^{n} \slash \forall q \in {\bf Z}, \;\;
\mid g^{q}(x) - x \mid^{2} > 0 \} \]

It is the open set formed by the elements which are not causally related to any 
$g$-iterates of themselves. We describe $\Omega_{g}$ in every conjugacy class:

\subsection{Elliptic case}

It is the case $gx = Rx + \tau$, where $R$ is a rotation fixing $\tau$.
Hence:

\[ \mid g^{q}x-x \mid^{2} = \mid (R^{q}-id)(x) + q\tau \mid^{2} =
\mid (R^{q}-id)(x) \mid^{2} + q^{2} \mid \tau \mid^{2} \]
since $\tau$ is orthogonal to $I$ which contains $\mid (R^{q}-id)(x) \mid^{2}$.
Remember also that $I$ is spacelike. Hence, we have three cases:

\begin{enumerate}

\item \emph{$\tau$ is timelike:\/} Then, $\Omega_{g} = \emptyset$
(Indeed, the term $q^{2}\mid \tau \mid^{2}$ is the leading term
since $(R^{q}-id)(x)$ is uniformly bounded).

\item \emph{$\tau$ is spacelike:\/} Then, $\Omega_{g} = {\bf M}^{n}$

\item \emph{$\tau$ is lightlike:\/} Then, 
the complement of $\Omega_{g}$ is the subspace formed by
$R$-periodic points, i.e., the $R^{n}$-fixed points.
This case contains a very special one, the linear case,
where $\tau$ is $0$. 

\end{enumerate}

\subsection{Loxodromic elements}
It is more suitable to write elements of \mink
in the form $(x,y,v)$ with $x$, $y$ in ${\mathbb R}$ and $v$ in ${\mathbb R}^{n-2}$,
so that the lorentzian form $Q$ is expressed by:

\[ Q = xy + \Vert v \Vert^{2} \]

where $\Vert \Vert^{2}$ denotes the usual euclidean norm on
${\mathbb R}^{n-2}$. Then, the action of of $L(g)$ is
defined by:

\[ L(g)(x,y,v) = (\lambda x, \lambda^{-1}y, Rv) \]

where $\lambda = e^{\zeta/2}$ is a positive real number, and $R$
a rotation of ${\mathbb R}^{n-2}$. Observe also that
$\tau$ can be assumed belonging to ${\mathbb R}^{n-2}$. Then:

 \begin{eqnarray*}
\mid g^{q}u - u \mid^{2} & = & q^{2} \mid\tau\mid^{2} + 
(\lambda^{q}-1)(\lambda^{-q}-1)xy + \Vert (R^{q}-id)(v) \Vert^{2} \\
 & = & q^{2} \mid\tau\mid^{2} - 4\mbox{sh}^{2}(q\zeta)xy + \Vert (R^{q}-id)(v) \Vert^{2}
 \end{eqnarray*}

The leading term $- 4\mbox{sh}^{2}(q\zeta)xy$ ensures:

\[ \Omega_{g} = \{ (x,y,v) \slash xy<0 \} \]

\subsection{Parabolic elements}

When $p$ is the nilpotent part of a unipotent element $id+p$ of
\son, we have:

\[ (id+p)^{k} = id + kp + \frac{k(k-1)}{2}p^{2} \]

When classifying parabolic elements, we proved that
the elliptic part satisfies $R \circ p = p$, and that the translation part
can be assumed fixed by $R$ and in $v_{1}^{\perp}$.
It follows:

\begin{eqnarray*} 
g^{q}(x) & = & R^{q}(id+p)^{q}(x) + q\tau + \sum_{k=0}^{q-1}kp(\tau) + 
\sum_{k=0}^{q-1}\frac{k(k-1)}{2}p^{2}(\tau) \\
 & = & R^{q}x + qp(x) + \frac{q(q-1)}{2}p^{2}(x) + q\tau + \frac{q(q-1)}{2}p(\tau) 
+ \frac{q(q-1)(q-2)}{6}p^{2}(\tau) 
\end{eqnarray*}

Hence, $x$ belongs to $\Omega_{g}$ if and only if for every $q$,
$T(q)$ is positive, where:

\[ T(q) = \mid \frac{R^{q}x-x}{q} + p(x) + \frac{(q-1)}{2}p^{2}(x) + \tau + 
\frac{(q-1)}{2}p(\tau) + \frac{(q-1)(q-2)}{6}p^{2}(\tau) \mid^{2} \]

{\bf The linear case $\tau=0$:\/}
the remaining terms all belong to $J^{\perp}$, and $p^{2}(x)$
belongs to $J$; moreover, $p(x)$ belongs to $I$ which is orthogonal
to $\frac{R^{q}x-x}{q}$ since this term belongs to $\mathcal I$.
Hence, $T(q)$ is equal to 
$\mid \frac{R^{q}x-x}{q} \mid^{2} + \mid p(x) \mid^{2}$ and is
therefore nonnegative. It vanish if and only if $\mid p(x) \mid^{2}=0$,
i.e., $x \in J^{\perp}$, and $\mid \frac{R^{q}x-x}{q} \mid^{2} =0$,
i.e. $x=R^{q}x$ since $\mathcal I$ is spacelike.

In other words, the complement of $\Omega_{g}$ is the set of
$R^{n}$-fixed points in the hyperplane $J^{\perp}$.

{\bf The tangent case $\tau \in J^{\perp} \setminus I$:\/}
Then, $p^{2}(\tau)=0$.
Observe that all the terms in $T(q)$ belong to $J^{\perp}$,
and that $p^{2}(x)$, $p(\tau)$ both belong to $J$. Thus, the
expression of $T(q)$ reduces to:

\[ \mid \frac{R^{q}x-x}{q} + p(x) + \tau\mid^{2} \]

The term $\frac{R^{q}x-x}{q}$ belongs to $\mathcal I$, and 
$p(x)+\tau$ belongs to ${\mathcal I}^{\perp}$. Therefore,
$T(q)$ is the sum of the nonnegative terms
$\mid\frac{R^{q}x-x}{q}\mid^{2}$ and $\mid\tau+p(x)\mid^{2}$. 
Moreover, since $\tau$ is assumed \emph{not} in $I$,
the second term never vanishes.

We conclude that in the tangent case, $\Omega_{g}$
is the entire \mink.

{\bf The transverse case $\tau \in {\bf M}^{n} \setminus J^{\perp}$:\/}
In this case, $p^{2}(\tau)$ is not $0$.
We develop $T(q)$, and seek for the leading terms:
since $p^{2}(\tau)$ is isotropic and orthogonal to
$p(\tau)$, $p(x)$, the leading terms are the terms in $q^{2}$:

\[ \frac{(q-1)^{2}}{4} \mid p(\tau)\mid^{2} + 
2\frac{(q-1)(q-2)}{6}\langle \tau \mid p^{2}(\tau) \rangle \]

According to $(\ast)$:

\[ \langle \tau \mid p^{2}(\tau) \rangle + \mid p(\tau)\mid^{2} + 
\langle p(\tau) \mid p^{2}(\tau) \rangle = 0\]

Hence, keeping in mind $\langle p(\tau) \mid p^{2}(\tau) \rangle = 0$,
the leading term is:

\[ \frac{3(q-1)^{2} - 4(q-1)(q-2)}{12}\mid p(\tau)\mid^{2} =
\frac{(q-1)(5-q)}{12}\mid p(\tau)\mid^{2} \]

Since $\mid p(\tau)\mid^{2} > 0$, we obtain that $T(q)$ is
negative for big $q$:

Conclusion: $\Omega_{g}$ is empty.

\subsection{Visibility of one point from another}
Points of the achronal domain studied above are 
points which cannot observe any $g$-iterate of themselves.
We wonder now, a base point $x_{0}$ being fixed,
how many $g$-iterates of a point $x$ are causally
related to $x_{0}$? For our purpose,
we just need to consider the case where $g$ is parabolic:

\begin{lema}
\label{paranongh}
Let $g$ be a parabolic isometry. Let $x_{0}$ be any point of \mink.
Then, there exist at least one point $x$ in \mink admitting
infinitely many $g$-iterates in the future or in
the past of $x_{0}$.
\end{lema}

\preua
We want to prove the existence of $x$ such that
for infinitely $q$, the norm $\mid g^{q}(x) - x_{0} \mid^{2}$
is negative. We have:

\[ \mid g^{q}(x) - x_{0} \mid^{2} = \mid x \mid^{2} + \mid x_{0} \mid^{2}
- 2 \langle g^{q}(x) \mid x_{0} \rangle \]

In the transverse case, the leading term in $q$ for 
$\langle g^{q}(x) \mid x_{0} \rangle = \langle g^{q}(x_{0}) \mid x \rangle$ is: 

\[ \langle \frac{q(q-1)(q-2)}{6}p^{2}(\tau) \mid x \rangle \]

When $x$ is not orthogonal to $p^{2}(\tau)$, this quantity
is positive with arbitrarly big values for infinitely
many $q$. Thus, 
such a $x$ admits infinitely $g$-iterates causally related
to $x_{0}$: the lemma is proven in this case.

In the linear or tangent cases, the leading term for 
$\langle g^{q}(x) \mid x_{0} \rangle$ is:

\[ \frac{q(q-1)}{2}\langle p^{2}(x)+p(\tau) \mid x_{0} \rangle \]

Therefore, the lemma is proven in this case by taking
$x$ for which $\langle p^{2}(x)+p(\tau) \mid x_{0} \rangle$
is positive. \fin

\section{Spacelike hypersurfaces and achronal domains}
\label{spapa}

From now, we consider a flat spacetime $M$ of dimension $n$,
and an isometric immersion $f: S \rightarrow M$ of a \emph{complete}
riemannian hypersurface. Let $\widetilde{M}$ be the universal covering
of $M$, ${\mathcal D}: \widetilde{M} \rightarrow {\bf M}^{n}$
the developping map, $\Gamma$ the 
fundamental group of $S$, and 
$\tilde{f}: \widetilde{S} \rightarrow \widetilde{M}$ a lifting
of $f$. There is an action of $\Gamma$ on $\widetilde{M}$, \emph{a priori}
non injective, such that $\tilde{f}$ is an equivariant map.

Let $\rho$ be the holonomy morphism of $M$ restricted to $\Gamma$:
it is a morphism $\rho: \Gamma \rightarrow \mbox{Isom}({\bf M}^{n})$.

We begin with an easy, but fundamental observation (see for example
\cite{mess}, or \cite{scannel}):

\begin{prop}
\label{project}
Let $\Delta$ be a timelike direction in ${\bf M}^{n}$.
Let $Q(\Delta)$ be the quotient linear space ${\bf M}^{n} \slash \Delta$,
and $\pi: {\bf M}^{n} \rightarrow Q(\Delta)$ the quotient map.
Equip $Q(\Delta)$ with the metric for which $\pi$ restricted to every
spacelike hyperplane orthogonal to $\Delta$ is an isometry.
Then, $\pi \circ {\mathcal D} \circ \tilde{f}$ is a homeomorphism.
\end{prop}

\preua
The map $\pi \circ {\mathcal D} \circ \tilde{f}$ is distance increasing,
and thus, a local homeomorphism. Moreover, it has the lifting
property: indeed, let $x$ be a point in $\widetilde{S}$, and 
$\bar{c}: [0,1] \rightarrow Q(\Delta)$ be a path such that
$\bar{c}(0) = \pi \circ {\mathcal D} \circ \tilde{f}(x)$.
Consider a maximal lifting $c: [0,t[ \rightarrow \widetilde{S}$:
for any sequence $t_{n}$ in $[0,t[$ converging to $t$,
the $c(t_{n})$ form a Cauchy sequence in $\widetilde{S}$
since $\pi \circ {\mathcal D} \circ \tilde{f}$ is distance increasing
and $\bar{c}(t_{n})$ is a convergent sequence. Therefore,
$c(t_{n})$ converge to some element $x_{\infty}$ of
$\widetilde{S}$. Therefore, $c$ can be continuously extended
by $c(t)=x_{\infty}$, which is a contradiction with the
definition of $t$, except if $t=1$.

Therefore, $\pi \circ {\mathcal D} \circ \tilde{f}$ is
a covering map. The proposition follows since $Q(\Delta)$
is simply connected.\fin

\begin{cora}
\label{quo}
The maps $f$ and ${\mathcal D} \circ \tilde{f}$ are embeddings. The natural morphism
$\Gamma \subset \pi_{1}(M)$ and the morphism
$\rho$ are injective. Moreover, 
${\mathcal D} \circ \tilde{f}(\widetilde{S})$ intersects every
timelike geodesic.\fin
\end{cora}

From now, we identify $\widetilde{S}$ with its image
under ${\mathcal D} \circ \tilde{f}$, and $\Gamma$ with its image
under $\rho$.

\begin{cora}
\label{quoiso}
Let $P_{0}$ be a lightlike affine hyperplane of ${\bf M}^{n}$,
with direction the orthogonal of a lightlike direction $\Delta_{0}$.
Then, either $\widetilde{S}$ does not intersect $P_{0}$, or
it intersects every affine line with direction $\Delta_{0}$
contained in $P_{0}$.
\end{cora}

\preua
The proof is completely similar to the proof of \ref{project};
the restriction to $S \cap P_{0}$ of the projection onto $P_{0}\slash\Delta_{0}$
is an isometry.\fin  

\begin{cora}
\label{abel}
If $\Gamma$ preserves a spacelike affine subspace ${\bf U}$, then
$\Gamma$ contains as a finite index subgroup a free abelian group of finite rank.
\end{cora}

\preua
Apply proposition \ref{project} with $\Delta \subset {\bf U}^{\perp}$.
The action of $\Gamma$ on $\pi({\bf U})$ has to be free and properly
discontinuous, and the same is true for the
action on ${\bf U}$. But this action
is isometric for the induced euclidean metric on ${\bf U}$;
the lemma then follows from Bieberbach's Theorem.\fin

As an immediate corollary of proposition \ref{project}, we obtain that
$\widetilde{S}$ is the graph of a contracting map 
$\varphi: {\mathbb R}^{n-1} \rightarrow {\mathbb R}$. Inversely,
it is worth to extend the definition of spacelike hypersurfaces
to (local) graphs of contracting maps from (euclidean) 
${\mathbb R}^{n-1}$ into (euclidean) ${\mathbb R}$.
In other words, it is not necessary to consider spacelike hypersurfaces
as smooth objects, Lipschitz regularity is enough. This low regularity
is sufficient to define a (non-riemannian) length metric on the spacelike
hypersurface, since Lipschitz maps are differentiable almost everywhere:
if $c: t \mapsto {\mathbb R}^{n-1}$ is a $C^{1}$-curve in ${\mathbb R}^{n-1}$,
define the length of $\varphi \circ c$ as the Lebesgue integral
of $\sqrt{1-\varphi'(c(t))^{2}} \Vert c'(t) \Vert$.

Another corollary is the following global property: 
spacelike hypersurfaces
are achronal closed set, meaning that if $x$ and $y$
are two distinct points in $\widetilde{S}$, one cannot
be contained in the past of future cone of the other, i.e., we have
$\mid x-y \mid^{2}>0$. Since $\widetilde{S}$ is $\gamma$-invariant,
it follows that it is contained in $\Omega_{\gamma}$
for every element $\gamma$ of $\Gamma$.
Actually, thanks to corollary \ref{quoiso}, we have a
quite better statement:

\begin{defina}
For every isometry $g$ of ${\bf M}^{n}$.
We define an open set $U(g)$ in the following way:

- $U(g) = \Omega_{g}$ if $g$ is loxodromic, or spacelike
elliptic, or non linear parabolic.

- $U(g) = \emptyset$ if $g$ is nonspacelike elliptic,

- $U(g) = {\bf M}^{n} \setminus J^{\perp}$ if $g$
is linear parabolic, where $J$ is the lightlike
affine line of fixed points of $g$,

\end{defina}

\begin{lema}
\label{SdanU}
The complete spacelike hypersurface $\widetilde{S}$ is contained in the interior of
$U(\Gamma) = \bigcap_{\gamma \in \Gamma} U(\gamma)$.
\end{lema}

\preua
First, we prove the inclusion $\widetilde{S} \subset U(\gamma)$
for every nontrivial element $\gamma$ of $\Gamma$. We observed
previously this inclusion when $U(\gamma) = \Omega_{\gamma}$;
thus, we have only $2$ cases to consider:

\begin{itemize}

\item \emph{the linear parabolic case:\/}
Since $\widetilde{S} \subset \Omega_{\gamma}$, $\widetilde{S}$
cannot intersect the lightlike line of fixed points $J$.
According to corollary \ref{quoiso}, it cannot intersect
$J^{\perp}$. 

\item \emph{the nonspacelike elliptic case:\/} 
In this case, $L(\gamma)$ admits a timelike line of fixed points.
According to corollary \ref{quo}, this line intersects $\widetilde{S}$
at some point $x$. Then, $\gamma x = x+\tau$: this is impossible
since $\tau$ is nonspacelike and that $\widetilde{S}$
is achronal.

\end{itemize}

Consider now a point $x$ in $\widetilde{S}$, and $B$ a small 
convex neighborhood of $x$ such that lightlike geodesic segment
contained in $B$ intersects $\widetilde{S}$. Then, for every
$\gamma$, the complement of $U(\gamma)$ is a union of lightlike
geodesic which cannot intersect $\widetilde{S}$, and thus, cannot intersect
$B$. It proves that $\widetilde{S}$ is contained in the interior
of $U(\Gamma)$.
\fin

\begin{defina}
The interior of $U(\Gamma)$ is denoted by $\Omega(\Gamma)$ and
called the achronal domain of $\Gamma$. The connected component
of the interior of $U(\Gamma)$ containing $\widetilde{S}$ is  
denoted by $\Omega(\widetilde{S})$.
\end{defina}

\rquea
\label{remplacefini}
For any isometry $g$, and for every nonzero integer $n$ we have $U(g^{n}) = U(g)$.
Therefore, if $\Gamma'$ is a finite index subgroup of $\Gamma$, the equality
$U(\Gamma') = U(\Gamma)$ holds. Thanks to this remark, we can replace 
$\Gamma$ by any finite index subgroup. 
\erquea

\rquea
\label{nini}
From now, we will assume that the group $\Gamma$ does not contain nonspacelike
elliptic and transverse parabolic elements, since we have shown that in
this case $\Gamma$ cannot be the holonomy group of a Cauchy-complete GH spacetime.
\erquea

\begin{prop}
\label{ghyp}
$\Omega(\widetilde{S})$ is a convex open set. It is
globally hyperbolic. The action of $\Gamma$ on 
it is free and achronal.
\end{prop}

\preua
The convexity follows from the convexity of connected components
of every $U(\gamma)$. The action of $\Gamma$ on $\Omega(\widetilde{S})$ 
is obviously free and achronal. The global hyperbolicity
follows from the following criteria (see \cite{beem}):
a spacetime is \emph{distinguished}
if for every pair $x \neq y$, the causal futures \emph{and} causal pasts
of $x$, $y$ are distinct.
Now, a distinguished spacetime is globally hyperbolic if for any 
pair of points $x$, $y$, the intersection
between the causal past of $x$ and the causal future of $y$
is empty or compact.

The open set $\Omega(\Gamma)$ is obviously distinguished,
we can thus apply the criteria above: let $x$, $y$
be two points in $\Omega(\Gamma)$, with $x$ in the causal future of
$y$. Then they are points $x'$, $y'$ near $x$, $y$ in $\Omega(\Gamma)$,
and such that $x$ is in the past of $x'$ and $y$ in the future
of $y'$ - we mean the causality relation in ${\bf M}^{n}$. 
Then, $F(y) \cap P(x)$ is contained in the interior
of $F(y') \cap P(x')$, and this interior is contained in
$\Omega(\Gamma)$ since $F(y') \cap P(x')$ is contained
in every $U(\gamma)$.\fin

From now, we assume that
$M$ is a Cauchy-complete globally hyperbolic spacetime admitting $S$ as a
complete Cauchy hypersurface of $M$.
The following observation is well-known:

\begin{prop}
$M$ is homeomorphic to $ S \times {\mathbb R}$, 
where every $S \times \{ \ast \}$ are Cauchy hypersurfaces. Any nonspacelike geodesic in $\widetilde{M}$
intersect $\widetilde{S}$. 
\end{prop}

\preua
See \cite{beem}. \fin

\begin{prop}
\label{danom}
The developing map ${\mathcal D}: \widetilde{M} \rightarrow {\bf M}^{n}$ is injective, with image contained in $\Omega(\widetilde{S})$.
\end{prop}

\preua
Fix  a timelike direction $\Delta$ in \mink. For any element $\tilde{x}$ of $\widetilde{M}$,
let $\delta(\tilde{x})$ be the line containing ${\mathcal D}(\tilde{x})$ with direction $\Delta$, and $d(\tilde{x})$ the connected component of ${\mathcal D}^{-1}(\delta(\tilde{x}))$ containing $\tilde{x}$.
Then, $d(\tilde{x})$ is a timelike geodesic in $\widetilde{M}$ and must therefore intersect
$\widetilde{S}$ at a single point $p(\tilde{x})$. 

If ${\mathcal D}(\tilde{x})= {\mathcal D}(\tilde{y})$, then $\delta(\tilde{x})=\delta(\tilde{y})=\delta$
and ${\mathcal D}(p(\tilde{x}))  = \delta \cap {\mathcal D}(\widetilde{S}) = {\mathcal D}(p(\tilde{y}))$. Hence, according to corollary \ref{quo}, $p(\tilde{x}) = p(\tilde{y})$,
and thus, $d(\tilde{x}) =d(\tilde{y})$
But the restriction of $\mathcal D$ to $d(\tilde{x})$ is an local homeomorphism from a topological line into ${\mathbb R}$: it is therefore injective. We obtain $\tilde{x}=\tilde{y}$,
i.e., $\mathcal D$ is injective.

The homeomorphism $M \approx S \times {\mathbb R}$ lift to a homeomorphism
$\widetilde{M} \approx \widetilde{S} \times {\mathbb R}$ 
where every $\widetilde{S} \times \{ t \}$ is a Cauchy hypersurface for $\widetilde{M}$.
Hence, for any $t$, every nonspacelike geodesic intersecting $\widetilde{S}$ meets
$\widetilde{S} \times \{ t \}$. It follows, thanks to proposition \ref{project}, that every nonspacelike line of \mink meets ${\mathcal D}(\widetilde{S} \times \{ t \})$. That's all we need to reproduce the arguments used in the preceding \S, leading to the conclusion
that for any $t$, ${\mathcal D}(\widetilde{S} \times \{ t \})$ is contained in
$\Omega(\Gamma)$. Since $\widetilde{M}$ is connected, 
we obtain that the image of $\mathcal D$ is contained in $\Omega(\widetilde{S})$. \fin

All the results above suggest a nice way towards the proof of
the main theorems: if the action of $\Gamma$ on
$\Omega(\widetilde{S})$ is proper, the quotient space 
$M(S) = \Gamma \backslash{\Omega}(\widetilde{S})$
is well-defined, and any globally 
hyperbolic manifold containing $S$ can
be isometrically embedded in $M(S)$. 
Unfortunately, things are not going so nicely: indeed, for
unipotent spacetimes, $\Omega(\widetilde{S})$ can be the entire Minkowski space
whereas the action of $\Gamma$ on \mink is not proper
(see remark \ref{unipcomplet}).

Anyway, this approach will be essentially successful after
ruling out a special case including unipotent spacetimes: 
the nonloxodromic case, i.e., the case where 
$\Gamma$ has no loxodromic element: this is the topic of
the next \S.

\rquea
\label{defdep}
Actually, when $M$ is maximal globally hyperbolic,
$\mathcal{C}(\widetilde{S})$, the
image of $\mathcal D$, is the so-called Cauchy domain, or domain of
dependence of the spacelike hypersurface $S$ as defined in
\cite{mess}, \cite{ande} or \cite{bonsante}: it is the
open set formed by points such that every lightlike geodesic through the
point intersects $S$. When the initial hypersurface $S$ is closed,
$\Omega(\widetilde{S})$ is the Cauchy domain
of $\widetilde{S}$ (see remark \ref{dependance}); but this is not true in the general case.
\erquea






\section{The nonloxodromic case }
\label{nol}
This is the case where $\Gamma$ has no loxodromic element.

\begin{prop}
\label{sanloxo}
If $\Gamma$ does not contain loxodromic elements, it is elementary;
i.e., up to finite index, its linear part is either contained in $SO(n-1)$ (up to conjugacy),
or its linear part stabilizes an isotropic vector.
\end{prop}

\preua
If $\Gamma$ has no loxodromic elements, according to \cite{benoist}, the same
is true for its Zariski closure $G$. Since $G$ has a finite number of connected
components, we can assume, restricting $\Gamma$ to a finite index subgroup
if necessary, that $\Gamma$ is contained in the identity component $G_{0}$.
If $G_{0}$ admits a parabolic element $g_{1}$ with fixed point $x$ in ${\bf S}^{n-2} = \partial{\bf H}^{n-1}$,
then $x$ is fixed by every parabolic element of $G_{0}$: indeed, if $g_{2}$ is another
parabolic element with fixed point $x_{2} \neq x_{1}$, for any small neighborhood $U$ of
$x_{2}$, $g_{1}^{k}(U)$ is arbitrarly near $x_{1}$, i.e. far from $x_{2}$ for big $k$.
In particular, $g_{1}^{k}(U)$ does not contain $x_{2}$. Therefore, for big $l$, $g_{2}^{l}g_{1}^{k}(U)$
is arbitrarly near $x_{2}$. In particular, if $k$ and $l$ are sufficiently big, the closure
of $g_{2}^{l}g_{1}^{k}(U)$ is contained in $U$, showing that $g_{2}^{l}g_{1}$ is loxodromic.

Hence, if $G_{0}$ admits a parabolic element, the fixed point of this parabolic
element is fixed by every element of $G_{0}$, meaning that $L(\Gamma)$ stabilizes a lightlike direction.

If not, $G_{0}$ admits only elliptic elements. Therefore, $1$-parameter subgroups
are all compact: it follows that $G_{0}$ is compact, hence, contained in a conjugate of $SO(n-1)$.\fin

According to Proposition
\ref{sanloxo}, the nonloxodromic case 
decomposes in the elliptic case, and the parabolic case:

\subsection{The elliptic case}

There is a flat euclidean metric on \mink preserved by $\Gamma$.
Remember also that according to \ref{SdanU},
elliptic elements of $\Gamma$ are all spacelike.
Hence, $U(\Gamma) = {\bf M}^{n}$, and $\Gamma$
acts freely on \mink. 
According to Bieberbach's Theorem, up to finite
index, there is a unique maximal affine subspace ${\bf U}$ on which $\Gamma$
acts by translations. Moreover, $SO(n-1) \subset SO(1,n-1)$ is the stabilizer of a point
in ${\bf H}^{n-1}$, i.e., preserves a timelike direction $\Delta_{0}$: it follows that 
${\bf U}$ is timelike. According to lemma \ref{suspendre}, 
$M(S)$ is thus a linear twisted
product over a translation spacetime.

\subsection{The parabolic case: }
\label{prab}
In this case, up to finite index, $L(\Gamma)$ preserves a isotropic vector $v_{0}$
and $\Gamma$ contains a parabolic element $\gamma_{1}$.
Denote by $P_{0}$ the orthogonal $v_{0}^{\perp}$, by $\Delta_{0} = P_{0}^{\perp}$
the direction containing $v_{0}$, and by 
$\overline{P}_{0}$ the quotient space $P_{0}\slash \Delta_{0}$. 
The
lorentzian quadratic form on \mink induces on $\overline{P}_{0}$
an euclidean quadratic form $\bar{q}_{0}$ which is preserved by the natural action
of $L(\Gamma)$.
The quotient space $Q_{0}={\bf M}^{n}\slash \Delta_{0}$ is a 
$(n-1)$-dimensional space foliated by affine planes 
with direction $\overline{P}_{0}$.  We denote
by $\pi_{0}: {\bf M}^{n} \rightarrow Q_{0}$ the quotient map.
Finally, we denote by $G(v_{0})$ the group of isometries of
\mink preserving $v_{0}$. 
The group $G(v_{0})$ has an induced affine action on $Q_{0}$:
there is a representation $\alpha: G(v_{0}) \rightarrow \mbox{Aff}(Q_{0})$ for which $\pi_{0}$ is
equivariant. The kernel of $\alpha$ is generated by translations along $\Delta_{0}$.

In a similar way, we consider the quotient map 
$\pi_{1}: {\bf M}^{n} \rightarrow {\bf M}^{n}\slash P_{0} = Q_{1}$:
the induced action of an element $g$ of $G(v_{0})$ on ${\bf M}^{n}\slash P_{0}$ is
a translation by a real number $b(g)$.
There is a canonical map $f: Q_{0} \rightarrow Q_{1}$.

We keep in mind that $\Gamma$ preserves the complete spacelike hypersurface
$\widetilde{S}$: this has some implications:

\begin{enumerate}

\item The projection $\pi_{0}(\widetilde{S})$ is an open set 
of the form $f^{-1}(I)$ where $I$ is a segment of the affine line $Q_{1}$ 
(cf. corollary \ref{quoiso}). 

\item $\Omega(\widetilde{S})$ is $\pi_{1}^{-1}(I')$ where $I'$ is a 
segment of $Q_{1}$
containing $I$. Indeed, every $U(\gamma)$
is a preimage by $\pi_{1}$ of a half-line in $Q_{1}$.
\item The action of $\Gamma$ on $\pi_{0}(\widetilde{S})$
\emph{via $\alpha$\/} is free and properly discontinuous (since the restriction 
of $\pi_{0}$
to $\widetilde{S}$ is a homeomorphism onto its image).


\end{enumerate}


Thus, point $(3)$ above implies that $\alpha(\Gamma)$ is a discrete subgroup
of $Aff(Q_{0})$. 

\subsubsection{A coordinate system on $Q_{0}$.}
Select a lightlike vector $v_{n}$ in \mink
transverse to $P_{0}$ and such that $\langle v_{0} \mid v_{n} \rangle = 1$,
and denote abusively $\overline{P}_{0}$ the intersection $P_{0} \cap v_{n}^{\perp}$.
Then, we decompose every element of \mink in $xv_{0}+z+yv_{n}$
with $z$ in $\overline{P}_{0}$.
The space $Q_{0}$ can be canonically decomposed in the form 
$\overline{P}_{0} \times {\mathbb R}$, so that the projection map $\pi_{0}$
is given by:

\[ \pi_{0}(xv_{0}+z+yv_{n}) = ( z , y) \]

Now, the action of the image under $\alpha$ of an element $g$ of $G(v_{0})$
acts on $Q_{0} \approx \overline{P}_{0} \times {\mathbb R}$ has the following expression:

\[ (\ast\ast) \;\;\;\;\;\;\;\;
\alpha(g)(z, y) = (R(g)(z)+u(g)+yv(g), y+b(g)) \]

where $u(g)$, $v(g)$ belong to $\overline{P}_{0}$, $R(g)$ is a rotation
in $\overline{P}_{0}$, and
$b(g)$ a real number. 

The action of the linear part of $g$ on 
${\bf M}^{n} \approx {\mathbb R} \oplus \overline{P}_{0} \oplus {\mathbb R}$
can be expressed by:

\[ g(xv_{0}+z+yv_{n}) = 
(x-\langle v(g) \mid R(g)(z) \rangle - \frac{y}{2}\mid v(g) \mid^{2})v_{0}+
(R(g)z+yv(g))+yv_{n} \]

the translation part being given by $\mu(g) v_{0} + u(g) + b(g)v_{n}$, where
$\mu(g)$ is not characterized by $\alpha(g)$. Observe that the
linear part $L(g)$ is uniquely defined by the induced isometry 
$\lambda(g)$ of $\overline{P}_{0}$,
i.e., by the pair $(R, v)$. An element $g$ is elliptic if $\lambda(g)$ admits a fixed point
in $\overline{P}_{0}$. If not, $g$ is parabolic.

The action on $Q_{0}$ of purely unipotent elements of $G(v_{0})$ has 
the following expression:

\[ (\ast\ast') \;\;\;\;\;\;\;\;
\alpha(g)(z, y) = (z+u(g)+yv(g), y+b(g)) \]

Therefore, they form a group ${\mathcal U}(v_{0})$.
Their action on \mink is:

\[ g(xv_{0}+z+yv_{n}) = 
((x-\langle v(g) \mid z \rangle - \frac{y}{2}\mid v(g) \mid^{2} + \mu(g))v_{0}+
(z+yv(g) + u(g))+(y+b(g))v_{n} \]

We distinguish the kernel ${\mathcal A}(v_{0})$ of 
$b: {\mathcal U}(v_{0}) \rightarrow {\mathbb R}$.
Every element $g$ of ${\mathcal A}(v_{0})$ is characterized by $u(g)$, $b(g)$,
$\mu(g)$. We recognize the group $\mathcal{A}$ discussed in the introduction for the definition
of unipotent spacetimes.

\rquea
\label{choixcoord}
The coordinate system is subordinated to the initial choice
of $v_{n}$ which is defined up to an element of $P_{0}$.
Let select $w$ an element of $P_{0}$, and consider the
new coordinates $(x', z', y')$, and the new morphisms
$u'$, $v'$, $\mu'$ defined with respect
to the decomposition 
$\Delta_{0} \oplus (P_{0} \cap (v_{n}+w)^{\perp}) \oplus \langle v_{n}+w \rangle$.
\erquea

\begin{itemize}

\item If $w$ is a multiple $av_{0}$, then 
$xv_{0}+z+yv_{n} = (x-ay)v_{0}+z+y(v_{n}+av_{0})$, showing that
$u'=u$, $v'=v$ and $\mu'=\mu$.

\item If $w$ belongs to $\overline{P}_{0} = P_{0} \cap v_{n}$,
then $y'=y$, $u'= u+ \langle u \mid w \rangle v_{0}$, 
$v' = v + \langle v \mid w \rangle v_{0}$. It follows:

\[ \mu'(g) = \mu(g) + \langle u(g) \mid w \rangle \]

\end{itemize}

\subsubsection{Triviality of $b$}
Assume the existence of some element $\gamma_{0}$
of $\Gamma$ for which $b(\gamma_{0}) \neq 0$.
Then, being invariant by a non-trivial translation, the
intervall $I = \pi_{1}(\widetilde{S})$ must be the whole
$Q_{1}$. According to corollary \ref{quoiso}, it follows
that $\widetilde{S}$ intersect every fiber of 
$\pi_{0}$. Remember now that we assume here the
existence of a parabolic element $\gamma_{1}$ in $\Gamma$,
and lemma \ref{paranongh}: select the base point
$x_{0}$ in $\widetilde{S}$ and consider the points
$x$ of \mink admitting infinitely many $\gamma_{1}$-iterates 
causally related to $x_{0}$. At one hand, according
to the proof of lemma \ref{paranongh}, suitable $x$
are the points for which $\langle p^{2}(x)+p(\tau) \mid x_{0} \rangle$
is positive, thus $x$ can be selected up to $\Delta_{0}$;
in particular, it can be selected in $\widetilde{S}$.
On the other hand, since they both belong to the $\Gamma$-invariant
hypersurface $\widetilde{S}$, no $\Gamma$-iterate of $x$ can
be causally related to $x_{0}$. We thus obtained
a contradiction.

\subsubsection{The purely unipotent case}
\label{prabb}
We have just proved that the morphism $b$ is trivial, i.e., 
the action of $\Gamma$ on $Q_{1}$ is trivial.
This case cannot arise when the initial spacelike hypersurface is compact.
Indeed, $\pi_{1}$ is then $\Gamma$-invariant, inducing a submersion
from $S$ onto ${\mathbb R}$. 

Modifying the coordinate system, we can assume that $0$ belongs to $I = \pi_{1}(\widetilde{S})$.
According to point $(3)$ above, $\Gamma$ acts freely and properly discontinuously on
$\overline{P}_{0} \times \{ 0 \}$ - observe that this action is the euclidean action
defined through $\lambda(\Gamma)$. 
According to Bieberbach's theorem, up to finite index, $\Gamma$ 
acts effectively by translations on some affine subspace
$\overline{\mathcal U}_{0} \times \{ 0 \} \subset \overline{P}_{0} \times \{ 0 \}$. In particular,
every $u(\gamma)$ belongs to the direction of $\overline{\mathcal U}_{0}$. 
But this observation works for any $t \in I$: there is some affine subspace
$\overline{\mathcal U}_{t} \times \{ t \} \subset \overline{P}_{0} \times \{ t \}$
on which the action of $\alpha(\gamma)$ has the expression:

\[ \alpha(\gamma)(z,t) = (z + u(\gamma) + tv(\gamma), t) \]

when $z$ belongs to $\overline{\mathcal{U}}_{t}$. It follows that the action
of every $\lambda(\gamma)$ on the linear space spanned by all the $\overline{\mathcal{U}}_{t}$
is unipotent; since $\overline{\mathcal{U}}_{0}$ is the \emph{maximal} $\lambda(\Gamma)$-unipotent
subspace, it must contain every $\overline{\mathcal{U}}_{t}$. Hence, every $v(\gamma)$
belongs to $\overline{\mathcal{U}}_{0}$: $\overline{\mathcal{U}}_{t}$ does not depend
on $t$. 

In the coordinate system ${\bf M}^{n} \approx {\mathbb R} \oplus \overline{P}_{0} \oplus {\mathbb R}$,
consider the subspace ${\mathcal U} = {\mathbb R} \oplus \overline{\mathcal{U}}_{0} \oplus {\mathbb R}$.
It is $\Gamma$-invariant, and the expression of the action of $L(\gamma) - id$ 
on ${\mathcal U}$ is:

\[ (L(\gamma) -id)(x,z,y) = 
(-\langle v(g) \mid z \rangle - \frac{y}{2}\mid v(g) \mid^{2}, z+yv(g), 0) \]

This action is unipotent. On the other hand, ${\mathcal U}$ is timelike.
Thanks to lemma \ref{suspendre}, and since the result we are proving
is up to linear twisted products, we can assume that 
${\mathcal U}$ is the entire Minkowski
space, i.e., $\Gamma$ is contained in the abelian group $\mathcal{A}(v_{0})$.

Since the action of $\Gamma$ on $\overline{P}_{0} \times \{ 0 \}$ is effective,
$u: \Gamma \rightarrow \overline{P}_{0}$ is injective.
We can thus parametrize $\Gamma$ by its $u$-translation vectors. Then,
if $E$ is the linear space spanned by the $u$-translation vectors, there is
a linear map $T$ from $E$ into $\overline{P}_{0}$ such that:

\[ \alpha(u)(z,y) = (z + u + yT(u),y) \]

Moreover, $\Gamma$ is  isomorphic to ${\bf Z}^{k}$. In particular, it is abelian:

\[ \mu(\gamma) + \mu(\gamma') - \langle u(\gamma') \mid v(\gamma) \rangle = \mu(\gamma\gamma') = \mu(\gamma'\gamma) = \mu(\gamma') + \mu(\gamma) - \langle u(\gamma) \mid v(\gamma') \rangle \]

It follows, for every $u$, $u'$ in $E$:
 
\[ \langle u \mid T(u') \rangle = \langle u' \mid T(u) \rangle \]

 The causality domain $U(u)$ is then
the open set formed by the points $(x, z, y)$ for which 
$u+yA(u) \neq 0$. The interval $I'$ (remember point $(2)$) is 
a connected component of the set 
$\{ y \in Q_{1} \approx {\mathbb R} \slash u+yA(u) \neq 0 \;\;\; 
\forall u \in \Gamma \}$

Denote by $y_{-}$, $y_{+}$ the extremities (maybe infinite)
of the intervall $I = \pi_{1}(\widetilde{S})$. Remember that we assume here
$y_{-} < 0 < y_{+}$. When $y_{\pm}$ are both finite, we change the coordinates
so that $y_{-} = -y_{+}$.

\begin{prop}
\label{limites}
At least one of the extremities $y_{\pm}$ is finite. If they are both finite then,
reducing as above to the case $y_{-} = -y_{+}$, for every $u$ in $E$ we have:

\[ \mid y_{+}T(u) \mid \leq \mid u \mid \]

If $y_{+} = +\infty$, then for every $u$ in $E$:

\[   \langle u \mid T(u) \rangle \geq -y_{-}\mid T(u) \mid^{2} \]

If $y_{-} = -\infty$, then for every $u$ in $E$:

\[   \langle u \mid T(u) \rangle \leq -y_{+}\mid T(u) \mid^{2} \]
\end{prop}

\rquea
\label{unipcomplet}
When every elements of $\Gamma$ are nonlinear unipotent isometries, 
$\Omega(\widetilde{S})$ is the entire Minkowski space. Then, according to
proposition \ref{limites}, $I$ and $I'$ are not equal.
\erquea

\preuda{limites}
Consider two elements $X = (x, y, z)$ and $X' = (x', y' z')$ of $\widetilde{S}$,
and any element $\gamma = g_{(u, v, \mu)}$ of $\Gamma$ (by definition, $v = T(u)$).
We want to evaluate the norm of $\gamma^{p}X-X'$: the first observation
is that $\mu(\gamma^{p}) = p \mu(\gamma)+\frac{p(1-p)}{2}\langle u \mid v \rangle$.
Therefore:

\[ \gamma^{p}(x,y,z) = (x+p \mu(\gamma)+\frac{p(1-p)}{2}\langle u \mid v \rangle-p\langle z \mid v\rangle - \frac{y}{2}p^{2} \mid v \mid^{2}, y, z+pu+ypv) \]

The leading term for $p$ increasing to infinity of the minkowski
norm of $\gamma^{p}X-X'$ is:
\[ p^{2}(\mid u \mid^{2} + (y+y')\langle u \mid v \rangle + yy'\mid v \mid^{2}) \]

Since no $\gamma$-iterate of $X$ can be in the causal future of $X'$, this term
must be positive for every $y$, $y'$ in $I$. Since at least one element of $\Gamma$
is parabolic, i.e., as a nonzero $v$-component, it follows immediatly that $I$ cannot
be the entire real line. When $y_{\pm}$ are both finite, $y_{-}=-y_{+}$, the limit
case $y=y_{+}$, $y'=y_{-}$ provides the required inequality for $u=u(\gamma)$.
This inequality extend to any $u$ in $E$ since $u(\Gamma)$ is a lattice of $E$.

When $y_{+} = +\infty$, the non-negativity of the leading term for increasing $y$
and $y'=y_{-}$ proves the required inequality for $u=u(\gamma)$. Once again,
the cocompactness of $u(\Gamma)$ in $E$ completes the proof.

Finally, the symetric case $y_{-}=-\infty$ admits a similar proof. \fin

\begin{prop}
\label{extension}
The map $T: E \rightarrow \overline{P}_{0}$ can be extended to a linear map
$\widehat{T}: \overline{P}_{0} \rightarrow \overline{P}_{0}$ satisfying the same properties, i.e.:

\begin{itemize}

\item $\overline{T}$ is symetric:
\[ \langle u \mid T(u') \rangle = \langle u' \mid T(u) \rangle \] 

\item If $y_{-} = -y_{+}$:

\[ \mid y_{+}\widehat{T}(u) \mid \leq \mid u \mid \]

\item If $y_{+} = +\infty$:

\[   \langle u \mid \widehat{T}(u) \rangle \geq -y_{-}\mid \widehat{T}(u) \mid^{2} \]

\item If $y_{-} = -\infty$:

\[   \langle u \mid \widehat{T}(u) \rangle \leq -y_{+}\mid \widehat{T}(u) \mid^{2} \]
\end{itemize}
\end{prop}

\preua
Consider first the case $y_{+}=+\infty$.
Let $F$ be the image of $T$, and $F^{\perp}$ the orthogonal of $F$. Let $K$ be the kernel of $T$, and $K^{\perp}$ the orthogonal of $K$
inside $E$. From the equalities $0 = \langle T(u) \mid u' \rangle = \langle u \mid T(u') \rangle$
for $u$ in $E$ and $u'$ in $K$, we see that $K \subset F^{\perp}$. 
On the other hand, for $u$ in $E \cap F^{\perp}$, we have 
$-y_{-}\mid T(v) \mid^{2} \leq \langle u \mid T(u) \rangle = 0$.
Therefore, $K = E \cap F^{\perp}$: the linear space $\overline{P}_{0}$ is the orthogonal sum $K^{\perp} \oplus F^{\perp}$. For $u$ in $K^{\perp}$ and $u'$ in $F^{\perp}$,
we define $\widehat{T}(u+u') = T(u)$. The linear map $\widehat{T}$ has the required properties.

The proof for the case $y_{-}=-\infty$ is completely similar. The last case to consider is
$y_{-}=-y_{+}$. We can assume without loss of generality that $y_{+}=1$.
We decompose $T(u)$ as the orthogonal sum $T_{0}(u) + B(u)$, where $T_{0}(u)$
belongs to $E$ and $B(u)$ to $E^{\perp}$. Let $B': E^{\perp} \rightarrow E$
the dual of $B$: it is uniquely defined by $\langle B'(v) \mid u \rangle = \langle v \mid B(u) \rangle$ for every $v$ in $E^{\perp}$ and every $u$ in $E$.

By hypothesis, $T_{0}$ is symetric and its sup-norm is less than one. Therefore, it is diagonalizable over ${\mathbb R}$, and its eigenvalues have all absolute values less than $1$.
For every symetric linear map $Z: E^{\perp} \rightarrow E^{\perp}$, define
$\widehat{T}_{Z}(u+v) = (T_{0}(u)+B'(v)) + (B(u)+Z(v))$. This is a symetric extension
of $T$, and we want to prove that $Z$ can be selected so that $\widehat{T}_{Z}$ has norm
less than $1$, i.e., that its eigenvalues have all absolute values less than $1$.

We first observe that the general case follows from the case where the sup-norm of $T_{0}$ is \emph{strictly} less than $1$. Indeed, the space of symetric maps with norm less than one is compact. Therefore, if for every positive integer $k$ the map
$(1-\frac{1}{k})T$ admits an extension $\widehat{T}_{k}$ with norm less than $1$,
we can extract from these $\widehat{T}_{k}$ a subsequence converging to
a symetric extension $\widehat{T}$ of $T$ which still has norm less than $1$.

Hence, we can assume that eigenvalues of $T_{0}$ have all absolute values strictly less than $1$. 
Let $\lambda$ be any real number with absolute value greater than $1$.
By hypothesis, $\lambda$ is not an eigenvalue of $T_{0}$, thus, we can define 
$Q(\lambda, v)=\lambda \mid v \mid^{2} - \langle (\lambda-T_{0})^{-1}B'v \mid B'v \rangle$.
For fixed $\lambda$, $Q(\lambda, v)$ is a quadratic form on $E^{\perp}$.
Therefore, $\langle Z(v) \mid v \rangle = \frac{Q(-1,v)+Q(1,v)}{2}$ defines a symetric linear map $Z$: our choice for extending $T$ will be $\widehat{T}=\widehat{T}_{Z}$.

We have the inequality $\mid T(u) \mid \leq \mid u \mid$ for every $u$ in $E$, i.e.:

\[\begin{array}{rcl}
\mid T_{0}(u)\mid^{2} + \mid B(u) \mid^{2} & \leq & \mid u \mid^{2} \\
\langle u \mid (1-T_{0}^{2})(u) \rangle & \geq & \mid B(u) \mid^{2}
\end{array} \]

Since $\pm 1$ is not an eigenvalue of $T_{0}$, we can apply the inequality above
to $u = (1-T_{0}^{2})^{-1}B'v$:

\[ \mid B(1-T_{0}^{2})^{-1}B'v \mid^{2} \leq \langle (1-T_{0}^{2})^{-1}B'v \mid B'(v) \rangle =
\langle B(1-T_{0}^{2})^{-1}B'v \mid v \rangle \leq \mid B(1-T_{0}^{2})^{-1}B'v\mid.\mid v \mid \]

Hence, for any $v$ in $E^{\perp}$:

\[ \mid B(1-T_{0}^{2})^{-1}B'v \mid \leq \mid v \mid \]

\[ \langle (1-T_{0}^{2})^{-1}B'v \mid B'v \rangle \leq \mid v \mid . \mid B(1-T_{0}^{2})^{-1}B'v\mid \leq \mid v \mid^{2} \]

It follows: 

\[ 2\mid v \mid^{2} \leq \langle [(1+T_{0})^{-1}+(1-T_{0})^{-1}] B'v \mid B'v \rangle \]

This last inequality is equivalent to $Q(-1, v)  \leq Q(1, v)$. Therefore, we have
$Q(-1, v)  \leq \langle Z(v) \mid v \rangle \leq Q(1, v)$. But, for a fixed $v$, $Q(\lambda, v)$
is an increasing function of $\lambda$ on $]-\infty, -1] \cup [1, +\infty[$. Hence, for any $\lambda$ in $]-\infty, -1[ \cup ]1, +\infty[$ and any nonzero $v$ in $E^{\perp}$,
we have the inequality:

\[  \langle Z(v) \mid v \rangle \neq Q(\lambda, v) \]

For $\lambda$ in $]-\infty, -1[ \cup ]1, +\infty[$, let $u+v$ be an element of the kernel of  $\lambda-\widehat{T}_{Z}$ with  $u$ in $E$ and  $v$
in $E^{\perp}$:

\[ \begin{array}{rcl}
T_{0}(u) + B'(v) & = & \lambda u \\
B(u) + Z(v) & =& \lambda v
\end{array} \]

Hence: 

\[ \begin{array}{rcl}
u & = & \lambda (\lambda-T_{0})^{-1}B'v  \\
Z(v) & =& \lambda v - B(\lambda-T_{0})^{-1}B'v 
\end{array} \]

Hence, $\langle Z(v) \mid v \rangle = Q(\lambda, v)$. According to the study above, it implies
$v=0$, and thus $u=0$: $\lambda$ is not an eigenvalue of $\widehat{T}_{Z}$.\fin

Since $\widehat{T}$ is symetric, there is an orthonormal basis $e_{1}, \ldots, e_{n-2}$ of $\overline{P}_{0}$ such that $T(e_{i}) = \lambda_{i}e_{i}$. For any collection
of real numbers $(\mu_{1}, \ldots, \mu_{n-2})$, we define a map
$\varphi: \overline{P}_{0} \approx {\mathbb R}^{n-2} \rightarrow {\mathcal A}(v_{0})$ by $\varphi(\sum t_{i}e_{i}) = 
g$ such that $u(g) = \sum t_{i}e_{i}$, $v(g) = \widehat{T}(u(g)) = \sum \lambda_{i}t_{i}e_{i}$ and $\mu(g) = \sum t_{i}\mu_{i} + \sum \lambda_{i}\frac{t_{i}(1-t_{i})}{2}$. It is easy to check that for any $(\mu_{1}, \ldots, \mu_{n-2})$, $\varphi$ is a  morphism.

 Let $\gamma_{1}, \gamma_{2}, \ldots, \gamma_{k}$ be generators of $\Gamma \approx {\bf Z}^{k}$. We can select $(\mu_{1}, \ldots, \mu_{n-2})$ so that for every $j$,
if $u(\gamma_{j}) = \sum t_{i}e_{i}$, then $v(\gamma_{j}) = T(u(\gamma_{j})) =
\sum \lambda_{i}t_{i}e_{i}$ and $\mu(\gamma_{j}) = \sum t_{i}\mu_{i} + \sum \lambda_{i}\frac{t_{i}(1-t_{i})}{2}$. Then, for every element $\gamma$ of $\Gamma$,
we have $\varphi(u(\gamma)) = \gamma$: the image of $\varphi$ is a Lie abelian subgroup
$A$ of ${\mathcal A}(v_{0})$ containing $\Gamma$
and isomorphic to $\overline{P}_{0} \approx {\mathbb R}^{n-2}$ considered as a group.

Observe that there is a vector $w$ in $\overline{P}_{0}$ such that
for every $i$, $\langle u_{ i}\mid w \rangle = -\mu_{i}-\frac{\lambda_{i}}{2}$. 
Therefore, according to remark \ref{choixcoord}, after a coordinate change,
we can actually assume $\mu_{i} = -\frac{\lambda_{i}}{2}$. In this coordinate system,
the $\mu$-component of $\varphi(\sum t_{i}e_{i})$ is 
$-\sum \frac{{\lambda}_{i}}{2}t_{i}^{2}$.

In summary, we have proven that $\Gamma$ is precisely as described in the introduction for the definition of unipotent spacetimes. Moreover, our choice of $\widehat{T}$ in proposition \ref{extension} makes sure that $\widetilde{S}$ is contained in
a connected component $\Omega$ of $\Omega(A)$ - cf. the introduction for the definition
of $\Omega(A)$.
But in the proof of proposition \ref{danom}, we observed that every lightlike geodesic
contained in $\widetilde{M}$ intersect $\widetilde{S}$: hence, 
${\mathcal D}(\widetilde{M})$ is contained in $\pi_{0}^{-1}(\pi_{0}(\widetilde{S}))$.
It follows that ${\mathcal D}(\widetilde{M})$ is contained in $\Omega$, i.e., that
$M$ can be embedded in the unipotent spacetime $\Gamma\backslash\Omega$.
This achieves the proof of Theorem \ref{caucomp} in the nonloxodromic 
parabolic case.

\section{The loxodromic case}
\label{loxocase}
It is the case where $\Gamma$ contains a loxodromic element $\gamma_{0}$.
We can then define the convex domain $\Omega_{lox}(\Gamma)$ which is 
the interior of the intersection between all the $U(\gamma)$ for \emph{loxodromic}
elements of $\Gamma$. 

\begin{lema}
\label{sanlox}
If $\Gamma$ contains loxodromic elements, 
$\Omega_{lox}(\Gamma) = \Omega(\Gamma)$.
\end{lema}

\preua
The inclusion $\Omega(\Gamma) \subset \Omega_{lox}(\Gamma)$ is obvious. 
The reverse inclusion would fail if for some nonloxodromic element $\gamma_{0}$
we had $\Omega_{lox}(\Gamma) \cap U(\gamma_{0}) \neq \Omega_{lox}(\Gamma)$.
Such a $\gamma_{0}$ cannot be tangent parabolic or spacelike elliptic since
in these cases $U(\gamma_{0})$ is the entire Minkowski space.
According to remark \ref{nini}, this $\gamma_{0}$ is actually linear
parabolic. Then, $U(\gamma_{0})$ is the
complement of some lightlike hyperplane $P(\gamma_{0})$. Thus,
since we assume $\Omega_{lox}(\Gamma) \cap U(\gamma_{0}) \neq \Omega_{lox}(\Gamma)$,
the convex domain $\Omega_{lox}(\Gamma)$ contains some $\gamma_{0}$-invariant
nonempty open subset $U$ of $P(\gamma_{0})$. 
Now, $P(\gamma_{0})$ contains a isotropic
direction $\Delta_{0}$ such that for every element
$x$ of $P(\gamma_{0}) \setminus \Delta_{0}$, the convex hull of the $\gamma_{0}$-orbit of
$x$ is a complete lightlike line $\Delta$ (with direction $\Delta_{0}$).
If $x$ is selected in $U \setminus \Delta_{0}$, $\Delta$ is contained in the convex domain
$U_{lox}(\Gamma)$. But this is impossible since the achronal domain $U(\gamma)$
of a loxodromic element cannot contain an entire lightlike affine line.
\fin

\begin{lema}
\label{loxregular}
$\Omega(\widetilde{S})$ is a regular convex domain.
\end{lema}

\preua
The achronal domain of a loxodromic element $\gamma$ admits two connected
component, one future complete, the other, past complete:
we denote them respectively by $U^{+}(\gamma)$, $U^{-}(\gamma)$.
Since $\Omega(\widetilde{S})$ is connected, it is contained in
only one of these components.
Up to time reversing isometries, we can assume that
there is some loxodromic element $\gamma_{0}$ for
which $\Omega(\widetilde{S}) \subset U^{+}(\gamma_{0})$. 

\emph{Claim:} for every loxodromic element $\gamma$, $\Omega(\widetilde{S})$
is contained in $U^{+}(\gamma)$.

Indeed, assume the existence of some
$\gamma_{1}$ for which $\Omega(\widetilde{S}) \subset U^{-}(\gamma_{1})$.
Then, $\Omega(\widetilde{S})$ must be contained in 
$U^{+}(\gamma_{0}) \cap U^{-}(\gamma_{1})$. But the reader can
easily check that in all cases, there is always a timelike
line avoiding the intersection $U^{+}(\gamma_{0}) \cap U^{-}(\gamma_{1})$.
This is a contradiction with Corollary \ref{quo}.

It follows that $\Omega(\widetilde{S})$ is the interior
of the intersection of all $U^{+}(\gamma)$ when
$\gamma$ describes all the loxodromic elements
of $\Gamma$. Denote by $\Lambda_{0}$ the set of repulsive
fixed points of loxodromic elements of $\Gamma$ in
${\mathcal J}$: we have just proved that 
$\Omega(\widetilde{S})$ is the interior of the future complete convex
domain defined by $\Lambda_{0}$. Denote now by
$\Lambda(\Gamma)$ the closure of $\Lambda_{0}$ in
${\mathcal J}$. Then, according to lemma
\ref{fermons}, $\Omega(\Lambda(\Gamma))$ contains
$\Omega(\widetilde{S})$: in particular, it is nonempty.
Then, according to lemma \ref{closecool}, $\Lambda(\Gamma)$
is future regular, i.e., $\Omega(\Lambda(\Gamma))$ is
a future complete regular convex domain.
The lemma follows then from lemma \ref{fermonsbis}\fin

\begin{cora}
The action of $\Gamma$ on $\Omega(\widetilde{S})$ is
free and properly discontinuous, and the
quotient space $M(S) = \Gamma \backslash \Omega(\widetilde{S})$
is a Cauchy-complete globally hyperbolic spacetime.
\end{cora}

\preua
According to propositions \ref{loxregular} and 
\ref{agissons}, the action of $\Gamma$
on $\Omega(\widetilde{S})$ is properly discontinuous. 
By definition of $\Omega(\widetilde{S})$, 
this action is free.
The corollary thus follows from proposition \ref{agissons}.
\fin

We may naively think that we achieved the proof of the
main theorem, but some point is missing: $M(S)$
is not a Cauchy-hyperbolic spacetime if the linear part
$\Gamma$ is not discrete in $SO_{0}(1,n-1)$.

\subsection{The elementary case:\/} 
\label{misnerr}
Here, we consider the case where, up to
finite index, $L(\Gamma)$ preserves some isotropic
direction $\Delta_{0}$. Since $\gamma_{0}$ is loxodromic, it admits a second
isotropic fixed direction: we denote this fixed direction by $\Delta_{1}$.
As in the nonloxodromic case, there are induced action of $\Gamma$
on ${Q}_{0}$, $Q_{1}$ the quotient spaces of \mink by $\Delta$ and $\Delta_{0}^{\perp}$.
The difference with the nonloxodromic case is that the induced
action on $Q_{0}$ is now given by:

\[ \alpha(\gamma)(z, y) = (R(\gamma)(z)+u(\gamma)+yv(\gamma), e^{a(\gamma)}y+b(y)) \]

The isometry $\gamma$ is loxodromic if and only if $a(\gamma) \neq 0$.

We choose as origin of $Q_{1}$ the fixed point of $\gamma_{0}$. Then, $b(\gamma_{0})=0$,
and $U(\gamma_{0})$ does not
contain the lightlike hyperplane $\pi_{1}^{-1}(0)$. Since $I = \pi_{1}(\widetilde{S})$, 
$I'=\pi_{1}(\Omega(\widetilde{S}))$ are
$\rho(\gamma_{0})$-invariant intervals of ${\mathbb R}$, 
they are both the interval 
$]-\infty, 0[$ or $]0, +\infty[$; let say, $I = I' = ]0, +\infty[$.
Moreover, $b(\gamma)=0$ for every $\gamma$ in $\Gamma$.

Let $\Gamma_{0}$ be the kernel of $a$: this is a discrete group of isometries of \mink.
For any element $\gamma'$ of $\Gamma_{0}$, and any element $\gamma$ of
$\Gamma$, the conjugate $\gamma\gamma'\gamma^{-1}$ is given  by:

\[ \alpha(\gamma\gamma'\gamma^{-1})(z,y) = (R(\gamma)R(\gamma')R(\gamma)^{-1}(z) + R(\gamma)(u(\gamma')) +
e^{-a(\gamma)}v(\gamma'), y) \]

Since $\alpha(\Gamma')$ is discrete, 
it follows that every $v(\gamma')=0$ since we can
choose as $\gamma$ a loxodromic element for which $a(\gamma)\neq 0$.
In other words, elements of $\Gamma_{0}$ are all elliptic.

The nondiscrete part $\Gamma^{nd}_{0}$
of $\Gamma_{0}$ has finite index in $\Gamma_{0}$ and is preserved
by conjugacy under $\Gamma$: the unique maximal $\Gamma^{nd}_{0}$-unipotent subspace ${\bf U}$ is
preserved by $\Gamma$. 

Now, we observe that $\Gamma_{0}$ preserves $\widetilde{S}$, has no loxodromic elements, and acts
trivially on $Q_{1}$: we recover the purely unipotent case. In particular, ${\bf U}$ is timelike.
Since it is $\Gamma$-invariant, up to a linear twisted product, 
we can assume thanks to
lemma \ref{suspendre} that
${\bf U}$ is the entire Minkowski space and that 
elements of $\Gamma^{nd}_{0}$ are spacelike translations. 

Consider now the action of $L(\Gamma)$ on the hyperbolic space ${\bf H}^{n-1}$:
since $L(\Gamma^{nd}_{0})$ is trivial, $L(\Gamma_{0})$ is a finite group.
Therefore, ${\mathcal H}$, the set of $L(\Gamma_{0})$-fixed points, is a $L(\Gamma)$-invariant
totally geodesic subspace of ${\bf H}^{n-1}$. It is not reduced to a point since
the isotropic direction $\Delta_{0}$ defines a $L(\Gamma_{0})$-fixed point in $\partial{\bf H}^{n-1}$.
Now, $L(\Gamma_{0})$ acts trivially on ${\mathcal H}$, and for any element $\gamma$ of $\Gamma$,
the commutator $\gamma\gamma_{0}\gamma^{-1}\gamma_{0}^{-1}$ is in the kernel of
$a$, i.e., in $\Gamma_{0}$. Therefore, $L(\gamma)$ is an isometry of the 
hyperbolic space 
${\mathcal H}$ commuting with the loxodromic element $L(\gamma_{0})$: 
hence, it preserves
the \emph{two} fixed points of $L(\gamma_{0})$ in $\partial{\mathcal H}$.
In other words, there is another  isotropic direction $\Delta_{1}$ 
preserved by the entire group $L(\Gamma)$.

Let $P$ be the spacelike direction $\Delta_{0}^{\perp} \cap \Delta_{1}^{\perp}$, and
$E$ the quotient space of \mink by $P$. Then, $\Gamma$ acts naturally on the $2$-plane
$E$, preserving the lorentzian metric induced from the metric on \mink. Denote by
$\pi_{P}: {\bf M}^{n} \rightarrow E$ the quotient map.
We parametrize $P$ by coordinates $x$, $y$ so that the induced lorentzian metric is
$dxdy$, and that the fixed $\gamma_{0}$-invariant lightlike hyperplanes are
$\pi_{P}^{-1}(x=0)$, $\pi_{P}^{-1}(y=0)$. Then, since the linear parts of loxodromic elements of
$\Gamma$ preserve the same isotropic directions than $L(\gamma_{0})$, the achronal
domain of any loxodromic element $\gamma$ of $\Gamma$ is $\pi_{P}^{-1}(\{ (x-x(\gamma))(y-y(\gamma)) < 0 \})$
for some real numbers $x(\gamma)$, $y(\gamma)$.
According to lemma \ref{sanlox}, it follows that $\Omega(\Gamma) = \Omega_{lox}(\gamma)$
has also the expression $\pi_{P}^{-1}(\{ (x-x_{0})(y-y_{0}) < 0 \})$.
Finally, since $\Omega(\Gamma)$ is $\gamma_{0}$-invariant, we have $x_{0} = y_{0} = x(\gamma) = y(\gamma) = 0$.

Then $\pi_{P}^{-1}(0)$ is a $\Gamma$-invariant spacelike subspace on which
$\Gamma$ acts by isometries. 
In other words, there is a coordinate system $(x, z, y)$ on \mink such that
the lorentzian quadratic form is $xy + \mid z \mid^{2}$ (where $\mid \mid$ is an euclidean norm)
and such that elements of $\Gamma$ acts according to the law:

\[ \gamma(x, z, y) = (e^{-a(\gamma)}x, R(\gamma)(z) + u(\gamma), e^{a(\gamma)}y) \]

As in corollary \ref{abel}, we see that the action on $\pi_{P}^{-1}(0)$ is free
and properly discontinuous. Once more, an appropriate use of Bieberbach's theorem
leads to the conclusion:
the quotient space $M(S) = \Gamma\backslash\Omega(\widetilde{S})$  
is a linear twisted product over a Misner spacetime.

\subsection{The non-elementary case}
\label{generic}
The last case to
consider is the case where no finite 
index subgroup of $L(\Gamma)$ fixes a isotropic direction.

Apply Theorems
\ref{ausl}, \ref{unip}: the non-discrete part 
$\Gamma_{nd}$ is nilpotent; and there is a unique maximal ${\bf U}$
on which the action of $\Gamma_{nd}$ is unipotent. 
Being unique, ${\bf U}$ is also
preserved by $\Gamma$. 

\begin{prop}
${\bf U}$ is timelike.
\end{prop}

\preua
If it is lightlike, then it contains a unique isotropic
direction which is $\Gamma$-invariant: contradiction.

If ${\bf U}$ is spacelike, then, according to corollary 
\ref{abel}, a finite index subgroup of $\Gamma$ is abelian. Now, an abelian subgroup
of $SO_{0}(1,n-1)$ fixes a pair of isotropic directions, except
if it consists only one elliptic elements fixing one and
only one timelike vector $v_{0}$. The first case is excluded by hypothesis,
and the second case is forbidden too, since in this case $v_{0}$ has to be
fixed by the entire group $\Gamma$, and, thus, to belong to ${\bf U}$, which
is a contradiction with the definition of ${\bf U}$.\fin

We denote by $L_{U}(\gamma)$
the linear part of the restriction of $\gamma$ to ${\bf U}$.

\begin{prop}
$L_{U}(\Gamma)$ is discrete.
\end{prop}

\preua
Since the action of $\Gamma$ on ${\bf M}^{n}\slash{\bf U}$ is a linear
euclidean action, the closure of $L_{U}(\Gamma)$ in $\mbox{SO}_{0}({\bf U})$
is the restriction to ${\bf U}$ of the closure of $L(\Gamma)$ in
$SO_{0}(1,n-1)$. Hence, the neutral component of this
closure is contained in $L_{U}(\Gamma_{nd})$, proving that this closure
is nilpotent and its elements 
does not have elliptic parts: it contains only hyperbolic or
unipotent elements. In particular, if it
is not trivial, its center contains a parabolic or hyperbolic
element. If it contains a parabolic element $\gamma_{0}$, the unique
isotropic direction of $\gamma_{0}$ is $L_{U}(\Gamma)$-invariant, and
thus $L(\Gamma)$-invariant. Contradiction.

If the center contains a hyperbolic element $\gamma_{0}$, then the 
two isotropic fixed points of $L(\gamma_{0})$
has to be exchanged by every element of $L(\Gamma)$. Contradiction.\fin

In other words, $L_{U}(\Gamma)$ is a non-elementary Kleinian group.
Actually, thanks to Proposition \ref{suspendre}, we can restrict our
study to the case ${\bf U} = {\bf M}^{n}$ (be aware that
$\Omega_{lox}(\Gamma)$ is obviously equal to
${\bf U}^{\perp} \oplus \Omega_{lox}(\Gamma_{\mid \bf U})$,
where $\Gamma_{\mid \bf U}$ denotes the restriction of
$\Gamma$ to ${\bf U}$). 

During the proof of lemma \ref{loxregular}, we actually proved
that $\Omega(\widetilde{S})$ is the regular convex domain
defined by the closure of repulsive points of loxodromic
elements of $\Gamma$. Therefore, it coincides precisely
with the regular convex domain $\Omega^{+}(\Lambda(\rho))$
appearing in definition \ref{definitioncauchyhyp}.

Hence, it seems that we have all elements in hand to
conclude, but one of them is missing:
$L$ might be noninjective! Let $\mathcal N$ be its
kernel: its elements are translations, and the translation vectors
form a lattice in some spacelike linear space $E$.
The action of $\Gamma$ on $\mathcal N$ by conjugacy
is linear; more precisely, the conjugacy by $\gamma$
maps the translation by vector $v$ to the translation
by vector $L(\gamma)v$. Hence, it extends to 
an isometric action on the euclidean space $E$,
which moreover preserves the lattice $\mathcal N$.
Hence, replacing $\Gamma$ by some finite index subgroup,
we can assume that this action is trivial.
Hence, 
$0 \rightarrow {\mathcal N} \rightarrow \Gamma \rightarrow L(\Gamma)
\rightarrow 0$ is a central extension.

Consider $p: {\bf M}^{n} \rightarrow Q$, 
the quotient of \mink by the spacelike subspace $E$:
$Q$ is naturally equipped with a Minkowski metric, and 
there is a $p$-equivariant action of $\Gamma$ on $Q$
which reduces to a (nonlinear) action of
$L(\Gamma) \approx \Gamma/{\mathcal N}$. 
$L(\Gamma)$ is still a discrete group
of isometries of $Q$, but now with injective linear part
morphism. 

On the other hand, since it is convex, and since
it is preserved by translations in $\mathcal N$, the regular
convex domain $\Omega(\widetilde{S})$ is preserved by translations
by vectors in $E$. 
Therefore, denoting $\Omega(L(\Gamma)) = p(\Omega(\widetilde{S}))$,
we have $\Omega(\widetilde{S}) = p^{-1}(\Omega(L(\Gamma)))$.
Finally, since $\Omega(\widetilde{S})$ is the regular convex
domain defined in \mink by the closure of repulisive fixed points
of $L(\Gamma)$, it should be clear to the reader that
$\Omega(L(\Gamma))$ is actually \emph{in the Minkowski space $E$\/}
the regular convex domain $\Omega(\Lambda(L(\Gamma)))$
as defined for the definition of Cauchy-hyperbolic spacetimes
(definition \ref{definitioncauchyhyp}). Hence, the
quotient $M(L) = L(\Gamma)\backslash\Omega(L(\Gamma))$ is
a Cauchy-hyperbolic spacetime. 

The quotient map $p$ induces now an isometric fibration 
$\bar{p}: M(S) \rightarrow M(L)$, with fibers isometric
to the torus $E/{\mathcal N}$.

\emph{Conclusion:\/} $M(S) = \Gamma\backslash\Omega(\widetilde{S})$
is a twisted product of a Cauchy-hyperbolic spacetime
by a flat torus.



\section{Summary of the proof}
\label{sumup}
The proof of theorem \ref{caucomp} is
quite intricate, and maybe not so easy to follow.
Thus, we consider usefull to summarize here
these proofs, and to add some comments.

\subsection{Proof of Theorem \ref{caucomp}}
We start with a Cauchy-complete GH spacetime $M$, with complete
Cauchy-surface $S$. In \S \ref{spapa}, we proved that the developping
map ${\mathcal D}: \widetilde{M} \rightarrow {\bf M}^{n}$
is injective, identifying $\widetilde{M}$ with an open 
domain contained in a open convex domain $\Omega(\widetilde{S})$,
this last domain being $\rho(\Gamma)$-invariant, where 
$\rho: \Gamma \rightarrow \mbox{Isom}({\bf M})^{n}$
is the holonomy morphism. Moreover, it is proved that
$\rho$ is injective: $\Gamma$ is then identified with
its image $\rho(\Gamma)$. These preliminaries would
conclude if $M(S) = \Gamma\backslash\Omega(\widetilde{S})$
was a model spacetime, but this is not true in general.

In \S \ref{nol}, we consider the case where $\Gamma$ 
has no loxodromic elements. We prove the folkloric
fact that the linear part $L(\Gamma)$ is then
elementary, i.e.,  preserves (up to finite index) a 
point in $\overline{\mathbb H}^{n-1}$ (this fact is maybe
most known when the group is discrete). Then, we prove 
that, if $M$ is not a linear twisted product over a translation
spacetime, 
then it is a linear twisted product over some spacetime
$M'$ for which the holonomy group is an abelian group acting
unipotently. All the difficulty is to extend this holonomy
group to an abelian Lie group $A$ of unipotent elements 
acting suitably on $\Omega(\widetilde{S})$
(Proposition \ref{extension}). It is then straightforward
to prove that $\Omega(\widetilde{S})$ is contained in a connected
component $\Omega$ of $\Omega(A)$. Therefore, $M'$ embedds isometrically
in the unipotent spacetime $\Gamma\backslash\Omega$.

The nonloxodromic case being ruled out, we assume then
the existence of a loxodromic element. We prove that the convex domain
$\Omega(\widetilde{S})$ is a regular convex domain, so that
the quotient $M(S)$ is indeed a spacetime, containing an embedded
copy of $M$. We consider once more a dichotomy: 

- either
$L(\Gamma)$ (up to index $2$) preserves a point 
in $\partial{\mathbb H}^{n-1}$: $M(S)$ is then a linear twisted product
over a Misner spacetime (\S \ref{misnerr}),

- either $L(\Gamma)$ is nonelementary: then, $M(S)$
is a linear twisted product over a spacetime satisfying the same properties,
with the additionnal requirement that the linear part of
the holonomy group is discrete. We then prove that it admits
a fibration by flat torii over a Cauchy-hyperbolic
spacetime (\S \ref{generic}).

\subsection{Proof of Theorem \ref{cauclo}}
Reconsider the proof above when $M$ is assumed Cauchy-compact.
A fundamental observation is that in a GH spacetime, every
closed spacelike hypersurface is necessarely a Cauchy surface,
hence, every GH spacetime containing an embedded copy of
$M$ is necessarely Cauchy-compact.

It follows that nontrivial linear twisted products 
destroy Cauchy-compactness, thus the ``up to linear twisted
products'' appearing in Theorem \ref{caucomp} can be erased
in the Cauchy-compact version \ref{cauclo}.
As observed previously, unipotent spacetimes cannot
be Cauchy-compact; thus, together with their finite coverings, they 
disappear in the Cauchy-compact version.

\rquea
\label{dependance}
In remark \ref{defdep}, we claimed that for maximal 
Cauchy-compact spacetimes, $\Omega(\widetilde{S})$
is the Cauchy domain of $\widetilde{S}$. This should be
obvious to the reader for the case of translation spacetimes,
and in the case of Misner spacetimes and twisted products
over standart spacetimes, it follows from \cite{bonsante}
and lemma \ref{loxregular}.
\erquea

\subsection{Absolute maximality}
\label{maxi}
Actually, proofs of theorems \ref{caucomp} and \ref{cauclo}
are not complete. Indeed, they also include the following statements:

(i) Translation spacetimes, Misner spacetimes, and twisted products of 
Cauchy-hyperbolic spacetimes by flat
torii are absolutely maximal, and any unipotent spacetime
can be tamely embedded in some absolutely maximal unipotent spacetime.

(ii) Any maximal Cauchy-compact GH spacetime is not only embedded
in a finite quotient of a translation spacetime, a Misner
spacetime or the twisted product of a standart spacetime by a 
flat torus, it is actually \emph{isometric} to one of them.

Actually, (ii) follows from (i), since, as noted previously,
closed spacelike hypersurfaces in GH spacetimes are always
Cauchy hypersurfaces. In other words, maximal Cauchy-compact
GH spacetimes are also absolutely maximal.

Let's prove (i). To avoid repetitions,
we call \emph{model spacetimes} the spacetimes listed in (i).
Let $M_{0}$ be one of them, and
consider an isometric embedding $f: M_{0} \rightarrow M$
in some GH spacetime $M$. Then, according to \ref{caucomp},
some finite covering $M'$ of $M$ tamely embeds in some
model spacetime $M_{1}$. Hence, some finite covering $M'_{0}$
of $M_{0}$ embeds isometrically in $M_{1}$.
By construction, every $M_{i}$ is the quotient of some
open domain convex $\Omega_{i}$ of \mink by the holonomy
group $\Gamma_{i}$. The embedding of $M'_{0}$ in
$M_{1}$ lifts to a isometric embedding 
$F: \Omega_{0} \rightarrow \Omega_{1}$. 
As any locally defined isometry between open subsets
of \mink, $F$ extends to a bijective isometry of
the entire Minkowski space: we thus can assume
that it is the identity map. Therefore,
$\Omega_{0} \subset \Omega_{1}$.
Moreover, the fundamental group
$\Gamma'_{0} \subset \Gamma_{0}$ of $M'_{0}$ is then
identified with a subgroup of $\Gamma_{1}$.

Now, remember that $F \;\;(=id)$ is the lifting of
the \emph{embedding} of $M'_{0}$ into $M_{1}$.
Hence, if the inclusion $\Omega_{0} \subset \Omega_{1}$ is
actually surjective, then $\Gamma'_{0} = \Gamma_{1}$,
and $f$ is a surjective embedding. 

Thus, our remaining task is to establish $\Omega_{0} = \Omega_{1}$,
except for the unipotent spacetimes. 

\emph{If $M_{0}$ is a translations spacetime:\/} 
then, $\Omega_{0}$ is the entire Minkowski space.
The equality $\Omega_{0} = \Omega_{1}$ follows.

\emph{If $M_{0}$ is loxodromic (Misner or twisted product of a 
Cauchy-hyperbolic spacetime by a flat torus):\/}
then, $\Gamma_{1}$ contains also loxodromic elements. 
By construction, $\Omega_{i}$ is the regular convex domain associated
to repulsive fixed points of loxodromic elements of $\Gamma_{i}$. 
Hence, $\Omega_{1} \subset \Omega_{0}$ follows from the inclusion
$\Gamma'_{0} \subset \Gamma_{1}$.

Let's now restrict our discussion on the remaining case, i.e.
the case where $M_{0}$ is unipotent. Then, since $\Omega_{1}$
contains $\Omega_{0}$ which is a domain in \mink between two
degenerate hyperplanes, or the half-space defined by a degenerate
hyperplane, it is obvious that $M_{1}$ is either a translation
spacetime, or an unipotent spacetime. But the first case
is excluded since $\Gamma_{1}$ contains parabolic elements (the
parabolic elements of $\Gamma_{0}$). Hence, $M_{1}$
is actually an unipotent spacetime. 

In section \ref{prabb}, we defined the following objects associated
to unipotent spacetimes: 

- the spacelike linear space $E_{i}$ generated by the
$u$-components of elements of $\Gamma_{i}$,

- a linear map $T_{i}: E_{i} \rightarrow {\mathbb R}^{n-2}$
expressing the $v$-components of elements of $\Gamma_{i}$
from their $u$-components.

Moreover, $\Omega_{i}$ has the form $\pi_{1}^{-1}(I_{i})$ where
$\pi_{1}$ is the projection map along $\Gamma_{i}$-invariant
degenerate hyperplanes from \mink onto a $1$-dimensional
linear space $Q_{i}$.

Consider first the case where
$I_{0}$ is a bounded interval $]y_{-}, y_{+}[$. As previously,
we assume without loss of generality the equality
$y_{+} = -y_{-}$. In proposition \ref{limites}, we proved that
the operator norm of $T_{0}$ is less than $\frac{1}{y_{+}}$.
Denote by $Y_{+}$ the inverse of this operator norm.
Then, proposition \ref{extension} implies that $T$ extends to some
symetric operator $\widehat{T}$, which provides some
abelian Lie group $A$ containing $\Gamma_{0}$ and
for which $\Omega  = \pi_{1}^{-1}(]-Y_{+}, Y_{+}[)$ is a connected
component of $\Omega(A)$.

Then, since $y_{+} \leq Y_{+}$, $M_{0}$ tamely embedds in 
the quotient $\Gamma_{0}\backslash\Omega$. On the other hand,
we claim that this last quotient is absolutely maximal.
Indeed, reconsider all the reasoning above, but now
assuming $M_{0} = \Gamma_{0}\backslash\Omega$, i.e.,
$y_{+} = Y_{+}$, or, equivalently, $\Omega_{0} = \Omega$.
Consider $\pi_{1}(\Omega_{1}) = ] y'_{-}, y'_{+}[$. Since
$\Gamma_{1}$ contains $\Gamma_{0}$, proposition \ref{limites}
implies $-Y_{+} \leq y'_{-}$ and $y'_{+} \leq Y_{+}$,
i.e., the equality  $\Omega_{0} = \Omega_{1}$ establishing
the absolute maximality of $M_{0}$.

The case where $I_{0}$ is not bounded can be treated
in a similar way.

\section{CMC foliations}
\label{CMC}

Let $M$ be a Cauchy-compact spacetime. A \emph{CMC foliation} on $M$
is a codimension $1$ foliation of $M$ for which all leaves
are spacelike and compact, 
with constant mean curvature. A \emph{CMC time function}
is a submersion $t: M \rightarrow {\mathbb R}$ such that:

- $t$ is increasing in time (i.e., its restriction
to any future oriented timelike curve is increasing)

- every fiber $t^{-1}(s)$ is a hypersurface with 
constant mean curvature, the CMC value being $s$.

Observe that in any spacetime admitting a CMC time function,
the fibers of this time function are the only CMC hypersurfaces
(see e.g. \cite{barzeg}). In particular, the CMC time function,
if it exists, is unique, and defines moreover the unique 
CMC foliation on the spacetime.

\begin{thma}
Every (flat) maximal Cauchy-compact GH spacetime admits 
a unique CMC foliation. It admits a CMC time function
if and only if it is finitely covered by a Misner
spacetime, or the twisted
product of a standart spacetime by a flat torus;
the CMC time function then takes value in $]-\infty, 0[$
(future complete case) or $]0, +\infty[$ (past complete case). 
Cauchy-compact translations spacetimes 
do not admit CMC time functions, and any CMC spacelike
closed hypersurface is a leaf of the unique CMC foliation.
\end{thma}

\preua
The case of standart spacetimes is treated in \cite{ande}.
It is also obvious that twisted products by flat torii
preserves the existence of CMC time functions.

In the two elementary cases (translation spacetimes and
Misner spacetimes) there is an abelian Lie group $A$ of
dimension $n-1$ acting by isometries, freely, properly discontinuously,
with closed spacelike orbits: these orbits are the
leaves of a CMC foliation. 

In the case of Misner spacetimes, the orbits of $A$
on the universal
covering are product of spacelike hyperbolae 
in ${\bf M}^{2}$ with euclidean spaces: 
there are thus obviously the fibers of some CMC time function.

In the case of translation
spacetimes, the CMC value of every
leaf is $0$ (they are actually totally geodesic).
Let $I$  be the quotient space of 
the action of $A$. For any closed spacelike CMC hypersurface
$S$, the projection of $S$ into  $I \approx \mathbb R$
is a compact intervall. It means that $S$ is tangent to
at least two $A$-orbits $O^{+}$, $O^{-}$, so that
$S$ is contained in the future of $O^{-}$ and in the past
of $O^{+}$. The max principle for CMC hypersurfaces
then implies that the constant mean value of $S$ is
bigger than the CMC value of $O^{+}$ and less than
the CMC value of $O^{-}$. Since $O^{\pm}$ are
both totally geodesic, the CMC value for $S$ is
actually $0$, i.e. $S$ is maximal.
Then, $S$ lifts as a maximal spacelike hypersurface 
$\widetilde{S}$ in
\mink: according to \cite{yau}, $\widetilde{S}$ is a parallel
spacelike hyperplane. In other words, $S$ is an
orbit of $A$.\fin

\medskip
\noindent
CNRS, UMPA, UMR 5669 \\
\'Ecole Normale Sup\'erieure de Lyon,\\
46, all\'ee d'Italie \\
69364 Lyon cedex 07, FRANCE \\

\noindent
LIFR M$I^{2}$P \\
Universit\'e Ind\'ependante de Moscou \\
Bol. Vlasievskii per. 11, Moscow, RUSSIA

\noindent
email: Thierry.Barbot@umpa.ens-lyon.fr \\
URL : http://umpa.ens-lyon.fr/

\end{document}